%% file: preprint.tex
\documentclass[a4paper]{article}

\usepackage[english]{babel}
\usepackage[utf8]{inputenc}
\usepackage{amsmath,amssymb,amsthm}
\usepackage{tensor}

\usepackage{url}

\usepackage[hidelinks]{hyperref}

\allowdisplaybreaks[3]

\theoremstyle{definition}
\newtheorem{dfn}{Definition}
\theoremstyle{plain}
\newtheorem{lem}[dfn]{Lemma}
\newtheorem{prop}[dfn]{Proposition}
\newtheorem{cor}[dfn]{Corollary}
\theoremstyle{remark}
\newtheorem*{rem}{Remark}
\newenvironment{prf}{\begin{proof}}{\end{proof}}

\input{commands}

\title{Killing spinor-valued forms and their integrability conditions}
\author{Petr Somberg, Petr Zima}
\date{}

\begin{document}

\maketitle

\begin{abstract}
\input{abstract}
\end{abstract}

\begin{quote}
\small
\textit{Keywords:}
Killing type equations, prolongation of differential systems,
projective invariance, spinor-valued differential forms, cone
construction, constant curvature space

\textit{MSC2010 classification:}
35N10, 53A20, 53A55, 53B20, 53C21, 58J70
\end{quote}

\input{body}

\paragraph{Acknowledgment:}

The authors gratefully acknowledge the support of the grants
GACR 306-33/1906357, GAUK 700217 and SVV-2017-260456.


\noindent
Petr Somberg, Petr Zima\\
Mathematical Institute of Charles University\\
Sokolovská 83, Praha 8 - Karlín, Czech Republic\\
E-mail: \url{somberg@karlin.mff.cuni.cz}, \url{zima@karlin.mff.cuni.cz}

\end{document}

%% file: commands.tex
\newcommand\er{\mathbb{R}}
\newcommand\ce{\mathbb{C}}
\newcommand\eps{\varepsilon}
\newcommand\ii{\mathsf{i}}
\newcommand\Tan{\mathrm{T}}
\newcommand\Vecf{\mathcal{X}}
\newcommand\Norm{\mathrm{N}}
\newcommand\V{\mathrm{V}}
\newcommand\W{\mathrm{W}}
\newcommand\K{\mathrm{K}}
\newcommand\Eb{\mathcal{E}}

\newcommand\Tra{\mathbb{T}}
\newcommand\Spnr{\Sigma}
\newcommand\Df{\Omega}
\newcommand\Sec{\Gamma}
\newcommand\Wedge{\textstyle\bigwedge}
\DeclareMathOperator*\End{End}
\DeclareMathOperator*\Cl{Cl}
\newcommand\intp{\mathbin\lrcorner}
\newcommand\clp{\cdot}
\newcommand\clf{{\gamma\clp}}
\newcommand\cclf{{\cone{\gamma}\clp}}
\newcommand\clv{{\gamma^\sharp\clp}}
\newcommand\cclv{{\cone{\gamma}^\sharp\clp}}
\newcommand\G{\mathrm{G}}
\newcommand\sola{\mathfrak{so}}
\newcommand\spinla{\mathfrak{spin}}
\newcommand\hola{\mathfrak{hol}}
\newcommand\covd{\nabla}
\newcommand\covdaff{\covd^\aff}
\newcommand\mcovd{\widehat\covd}
\newcommand\mcovdaff{\mcovd^\aff}
\newcommand\covdv{\covd^\V}
\newcommand\kcovd{\widetilde\covd}
\newcommand\covdg{\covd^g}
\newcommand\covda{\covd^a}
\newcommand\covdap{\covd^{a'}}
\newcommand\covdtra{\covd^\Tra}
\newcommand\covdtrav{\covd^{\Tra\V}}
\newcommand\dif{\mathrm{d}}
\newcommand\difV{\dif^\V}
\newcommand\difg{\dif^g}
\newcommand\difa{\dif^a}
\newcommand\dr{\dif r}
\newcommand\vr{\partial_r}
\newcommand\curv{\mathcal{R}}
\newcommand\curvaff{\curv^\aff}
\newcommand\curvv{\curv^\V}
\newcommand\mcurv{\widehat\curv}
\newcommand\kcurv{\widetilde\curv}
\newcommand\curvg{\curv^g}
\newcommand\curva{\curv^a}
\newcommand\weyl{\mathcal{W}}
\newcommand\Rho{\mathcal{P}}
\newcommand\cone{\overline}

\newcommand\isoto{\xrightarrow{\sim}}
\DeclareMathOperator*\tsum{\textstyle\sum}
\DeclareMathOperator\rank{rank}
\DeclareMathOperator\Scal{Scal}
\DeclareMathOperator\Shape{S}
\newcommand\aff{\mathsf{aff}}
\newcommand\Tw{{\frac{3}{2}}}
\newcommand\nh{\nobreakdash-}
\newcommand\Me{M_\eps}

%% file: abstract.tex
We study invariant systems of PDEs defining Killing vector-valued
forms, and then we specialize to Killing spinor-valued forms.
We give a detailed treatment of their prolongation and integrability
conditions by relating the point-wise values of solutions to the curvature
of the underlying manifold.
As an example, we completely solve the equations on model spaces of
constant curvature producing brand new solutions which do not come
from the tensor product of Killing spinors and Killing-Yano forms.

%% file: body.tex
\section{Introduction}

Killing equations are a class of invariant over-determined systems
of partial differential equations, appearing naturally in many problems
related to (pseudo\nh) Riemannian geometry.
One of the most prominent examples are the Killing vectors,
corresponding to infinitesimal isometries of Riemannian manifolds.
In the present article we focus on another specific example in the
hierarchy of Killing equations, termed \emph{Killing spinor-valued forms}.
We introduce relevant Killing equations and deduce their properties
mostly implied by integrability of the differential system in question.
We shall start our analysis in a rather general context, and then gradually
specialize to the cases of most author's interest.
As an application of general results we shall completely resolve the
Killing equations on model spaces of constant curvature.

The main motivation for the study of Killing spinor-valued forms is that
they are a natural generalization of both Killing spinors and
Killing-Yano forms.
The Killing spinors and Killing-Yano forms play a dominant role in the
geometrical analysis on Riemannian manifolds, e.g.\ the study of Dirac
and Laplace operators and the associated eigenvalue problems.
Subsequently, the two examples of Killing type equations gained their
own interest in theoretical physics, too.

A central question in the subject asks for (pseudo\nh) Riemannian
manifolds admitting non-trivial solutions of Killing type equations,
and their relation to the underlying geometric structure for which
they occur.
To some extent, this question is answered by the \emph{integrability
conditions} which relate the solutions with the curvature properties of
manifolds.
Moreover, the Killing spinors and Killing-Yano forms are closely
related to special Riemannian structures, e.g.\ Sasakian and
$\G_2$-manifolds.

The general interest in Killing spinors was stimulated and
accelerated by T.\ Friedrich's inequality \cite{fri80} for
the eigenvalues of Dirac operator.
He also proved that a Riemannian manifold admitting Killing spinors is
Einstein, which is a direct consequence of the first integrability
condition.
The so called \emph{cone construction}, cf.\ \cite{bar96}, relates
Killing spinors with parallel spinors on the metric cone, and thus
allows a classification via holonomy of the cone.

Killing-Yano forms were introduced by Yano in \cite{yan52} within the
framework of his study of Killing vectors and harmonic tensors.
Integrability conditions and cone construction for (special)
Killing-Yano forms were deduced in \cite{sem03}.
The relevant Killing equations are also examples of the so called
\emph{first BGG operator} in projective geometry.
Namely, they can be efficiently described in the context of parabolic
geometries using tractor calculus, cf.\ \cite{hsss12a}, \cite{hsss12b}.

Killing spinor-valued forms already appeared in theoretical physics in
the construction of Kaluza-Klein supergravity, cf.\ \cite{dnp86},
\cite{dp83}.
A systematic treatment of Killing spinor-valued forms can be found in
\cite{sz16}, the main result being the cone construction for special
Killing spinor-valued forms.
The present paper is a continuation of this effort.

Here is a brief summary of the content of our article.
In the next Section~\ref{sec:kvf}, we start by deducing a
\emph{prolongation} of the defining Killing equation.
In general, the prolongation procedure transforms the original
differential system into a closed one by introducing new
indeterminate variables for undetermined components of
first derivatives.
This corresponds to certain extension of the initial bundle to a larger one
equipped with suitable connection, such that the original system of
equations is equivalent to the equation for parallel sections with
respect to the newly constructed connection.
The prolongation allows to write down the integrability conditions in an
explicit way.
Our approach is based on direct computations that are guided by
representation theoretical considerations.
We shall also analyze intrinsic projective invariance of Killing-Yano
forms, and compare our results with those obtained by more abstract
methods based on the tractor~calculus.

We shall consider vector-valued differential forms that
take values in an arbitrary vector bundle equipped with linear
connection.
The presence of vector values, when compared to the scalar-valued
case, yields additional terms induced by the curvature acting just
on the values.
Later on we shall specialize to (pseudo\nh) Riemannian manifolds and
spinor-valued differential forms.
Note that all our results are valid in arbitrary signature of the
metric.
On the other hand, we will not attempt to cover other generalizations
such as affine connections with torsion or conformal Killing equations
that would complicate our computations.

In Section \ref{sec:skvf} we shall generalize \emph{special}
Killing-Yano forms that are accommodated to the cone construction
equivalence mentioned above.
An example of this kind is the contact form on a Sasakian
manifold.

Section \ref{sec:ksf} is devoted to considerations of spinor-valued
differential forms and specializing former achievements into the
language of spinor calculus.

In Section \ref{sec:cone} we briefly review the cone construction and
the already known equivalences for Killing equations.
Namely, we present explicit formulas for solutions which we later use in
the example of spaces of constant curvature.

In the last Section \ref{sec:ccurv} we employ the cone construction and
the integrability conditions with the aim to describe all Killing
spinor-valued forms on model spaces of constant curvature.
We discover additional solutions in degree $1$ and describe
them explicitly by means of a new variant of the cone construction.
These solutions are the first example of Killing spinor-valued forms
that are not spanned by tensor product of Killing spinors and
Killing-Yano forms.

\section{Killing vector-valued differential forms}
\label{sec:kvf}

Let $(M, \covdaff)$ be a smooth manifold of dimension~$n$, $\covdaff$
a torsion-free affine connection and $(\V, \covdv)$ a real or complex
vector bundle over $M$ equipped with a linear connection. We denote by
$\covd$ the linear connection combined from $\covdaff$ and $\covdv$,
acting on mixed tensors built out of the tangent bundle $\Tan M$ and $\V$
by means of duals and tensor products.
The situation of most interest for us are the vector-valued, or more
specifically, $\V$\nh valued differential forms which we will call just
$\V$\nh valued forms for short.
We use the notation based on superscripts
indicating the origin of all objects involved, e.g.\ the curvature
operators $\curvaff$, $\curvv$ and $\curv$ associated to $\covdaff$,
$\covdv$ and $\covd$, respectively.

\begin{dfn} \label{kvf}
Let $\Phi$ be a $\V$\nh valued form of degree $p \in \{0, \dots, n\}$.
Then $\Phi$ is a \emph{Killing $\V$\nh valued form} provided there
exists a $\V$\nh valued form $\Xi$ of degree $p+1$ such that
\begin{align} \label{eq:kvf}
  \covd_X \Phi &= X \intp \Xi,
  & & \text{for all } X \in \Vecf(M),
\end{align}
where $X \intp \Xi$ denotes the contraction of the form $\Xi$ by vector
field $X$.
\end{dfn}

In other words, $\Phi$ is Killing if and only if its covariant
derivative is totally skew-symmetric. The form $\Xi$ is hence uniquely
determined as the normalized skew-symmetrization of the covariant
derivative,
\begin{align} \label{eq:kvfskew}
  \Xi &= \tfrac{1}{p+1}\, \text{skew-symm.} (\covd \Phi)
    = \tfrac{1}{p+1}\, \difV \Phi,
\end{align}
which equals the exterior covariant derivative of\/ $\Phi$.
Because the equation \eqref{eq:kvf} does not imply for $p = 0$ any
restriction on $\Phi$, we will assume $p\ge 1$ for the rest of this
section.
On the other hand, for $p = n$ we set $\Wedge^{n+1} \Tan^*M$ to be the
zero vector bundle and \eqref{eq:kvf} is thus equivalent to $\Phi$ being
parallel.

In the scalar-valued case the solutions are often termed
\emph{Killing-Yano forms}, and we stick to this terminology
in order to clearly distinguish between the
scalar-valued case and the general vector-valued one.
It is well known that \eqref{eq:kvf} is a
projectively invariant system of partial differential
equations, and hence an invariant of the projective class of affine
connections $[\covdaff]$. A more detailed discussion around this observation
is given in Section \ref{sec:proj}.

\subsection{Killing connection}

In order to deduce a prolongation and integrability conditions for
the differential system \eqref{eq:kvf}, we decompose the action of curvature on $\Phi$
into components according to their tensor symmetry types.
This approach yields the prolongation in an invariant form, see Section
\ref{sec:proj} for further discussion.
In fact, the value in the vector bundle $\V$ does not play a serious
role and the whole procedure is parallel to the one for the
scalar-valued Killing-Yano forms, cf.\ \cite{sem03}.
We have
\begin{align} \label{eq:tfcurv}
  \curv_{X,Y} \Phi &= \covd^2_{X,Y} \Phi - \covd^2_{Y,X} \Phi ,
\end{align}
which is a $\V$\nh valued covariant $(p+2)$\nh tensor skew-symmetric
separately in the first two and the remaining $p$ indices.
By the Littlewood-Richardson rule the corresponding decomposition is
\begin{align} \label{eq:decomp2}
  (2) \otimes (p) &\simeq (p,2) \oplus (p+1,1) \oplus (p+2).
\end{align}
This result corresponds to the tensor product decomposition of
irreducible representations for the general linear group.
Here we denote by $(c_1,c_2,\dots)$ the space of tensors with symmetries
corresponding to the \emph{dual} partition, i.e.\ the Young diagram with
$c_1$ boxes in the first column, $c_2$ boxes in the second column, etc.
For example, $(p+2)$ is the totally skew-symmetric component, on the
other hand, $(p,2)$ is the component such that the skew-symmetrization
over any subset of $p+1$ indices vanishes. In the case $p = 1$ the
decomposition degenerates and the term $(p,2)$ disappears.

Now suppose that $\Phi$ is a Killing $\V$\nh valued form. From
\eqref{eq:kvf} and \eqref{eq:tfcurv} we get
\begin{align} \label{eq:kvfcurv}
  \curv_{X,Y} \Phi
    &= Y \intp (\covd_X \Xi) - X \intp (\covd_Y \Xi).
\end{align}
The right hand side depends linearly just on the first covariant
derivative of\/ $\Xi$, which in general decomposes according to
\begin{align} \label{eq:decomp1}
  (1) \otimes (p+1) &\simeq (p+1,1) \oplus (p+2).
\end{align}
Comparing the decompositions \eqref{eq:decomp2} and \eqref{eq:decomp1}
we conclude that the $(p,2)$\nh type component of\/ $\curv \Phi$
vanishes, and $\curv \Phi$ may be sufficient for computing the covariant
derivative of\/ $\Xi$.
In what follows we confirm these ideas and deduce explicit~formulas.

Let us denote the partially and totally skew-symmetrized action of the
curvature on skew-symmetric forms
\begin{align} \label{eq:skewcurv1}
  \curv_X \wedge \Phi
    &= \tsum_{j=1}^n e^j \wedge (\curv_{X,e_j} \Phi),
\\\label{eq:skewcurv2}
  \curv \wedge \Phi
    &= \tfrac{1}{2}
        \tsum_{i=1}^n e^i \wedge (\curv_{e_i} \wedge \Phi)
    = \tfrac{1}{2}
        \tsum_{i,j=1}^n e^i \wedge e^j \wedge (\curv_{e_i,e_j} \Phi),
\end{align}
where $\{e_1,\dots,e_n\}$ is a tangent frame and
$\{e^1,\dots,e^n\}$ its dual coframe. The three components with
different tensor symmetries are given by
\begin{align} \label{eq:curvdecomp1}
  (\curv \Phi)^{(p+2)}_{X,Y}
    &= \tfrac{2}{(p+1)(p+2)}\, Y \intp (X \intp (\curv \wedge \Phi)),
\\\label{eq:curvdecomp2}
\begin{split}
  (\curv \Phi)^{(p+1,1)}_{X,Y}
    &= \tfrac{1}{p} \Big( Y \intp (\curv_X \wedge \Phi)
          - X \intp (\curv_Y \wedge \Phi) - {}
\\
    &\qquad
          - \tfrac{4}{p+2}\, Y \intp (X \intp (\curv \wedge \Phi))
        \Big),
\end{split}
\\\label{eq:curvdecomp3}
\begin{split}
  (\curv \Phi)^{(p,2)}_{X,Y}
    &= \curv_{X,Y} \Phi - (\curv \Phi)^{(p+1,1)}_{X,Y}
      - (\curv \Phi)^{(p+2)}_{X,Y} {}
\\
    &= \curv_{X,Y} \Phi
      - \tfrac{1}{p} \Big( Y \intp (\curv_X \wedge \Phi)
          - X \intp (\curv_Y \wedge \Phi) - {}
\\
    &\qquad
          - \tfrac{2}{p+1}\, Y \intp (X \intp (\curv \wedge \Phi))
        \Big).
\end{split}
\end{align}
A straightforward computation verifies that the components indeed have
the appropriate symmetries. It
also easily follows that the $(p,2)$\nh type component vanishes
automatically for $p = 1$.

\begin{prop} \label{kvfprol}
Let $\Phi$ be a Killing $\V$\nh valued form of degree $p \ge 1$ and $\Xi$
the corresponding $\V$\nh valued form of degree $p+1$. Then it holds
\begin{align} \label{eq:kvfprol}
  \covd_X \Xi
    &= \tfrac{1}{p} \Big(
        \curv_X \wedge \Phi
        - \tfrac{1}{p+1}\, X \intp (\curv \wedge \Phi)
      \Big),
  & & \text{for all } X \in \Vecf(M),
\end{align}
as well as
\begin{align} \label{eq:kvfint0}
  (\curv \Phi)^{(p,2)}_{X,Y} &= 0,
  & & \text{for all } X, Y \in \Vecf(M).
\end{align}
Equivalently, $\curv_{X,Y} \Phi$ is completely determined by
$\curv_X \wedge \Phi$.
\end{prop}

\begin{prf}
By \eqref{eq:kvfcurv}, \eqref{eq:skewcurv1}, \eqref{eq:skewcurv2}, we have
\begin{align*}
\begin{split}
  \curv_X \wedge \Phi
    & = \tsum_{j=1}^n e^j \wedge (e_j \intp (\covd_X \Xi)
        - X \intp (\covd_{e_j} \Xi))
\\
    & = (p+1) \covd_X \Xi
      - \tsum_{j=1}^n (\langle e^j, X \rangle \covd_{e_j} \Xi
        - X \intp (e^j \wedge (\covd_{e_j} \Xi)))
\\
    & = p \covd_X \Xi
      + X \intp \Big(
          \tsum_{j=1}^n e^j \wedge (\covd_{e_j} \Xi)
        \Big),
\end{split}
\\
\begin{split}
  \curv \wedge \Phi
    & = \tfrac{1}{2} \tsum_{i=1}^n e^i \wedge \Big( p \covd_{e_i} \Xi
        + e_i \intp \Big( \tsum_{j=1}^n e^j \wedge (\covd_{e_j} \Xi)
      \Big) \Big)
\\
    & = (p+1) \tsum_{j=1}^n e^j \wedge (\covd_{e_j} \Xi),
\end{split}
\end{align*}
proving \eqref{eq:kvfprol}. Then we substitute \eqref{eq:kvfprol}
into \eqref{eq:kvfcurv},
\begin{align*}
\begin{split}
  \curv_{X,Y} \Phi
    & = Y \intp \Big( \tfrac{1}{p} \Big( \curv_X \wedge \Phi
        - \tfrac{1}{p+1}\, X \intp (\curv \wedge \Phi) \Big) \Big) - {}
\\
  & \quad
      - X \intp \Big( \tfrac{1}{p} \Big( \curv_Y \wedge \Phi
        - \tfrac{1}{p+1}\, Y \intp (\curv \wedge \Phi) \Big) \Big)
\\
    & = \tfrac{1}{p} \Big(
      Y \intp (\curv_X \wedge \Phi) - X \intp (\curv_Y \wedge \Phi)
      - \tfrac{2}{p+1}\, Y \intp (X \intp (\curv \wedge \Phi)) \Big),
\end{split}
\end{align*}
proving \eqref{eq:kvfint0}. The last statement is a consequence of
\eqref{eq:kvfint0} and the fact that $\curv \wedge \Phi$ is just
skew-symmetrization of\/ $\curv_X \wedge \Phi$.
\end{prf}

The prolongation of vector-valued Killing forms then easily follows from
\eqref{eq:kvfprol}.
The appropriate prolongation vector bundle is the direct sum of $\V$\nh valued $p$\nh form
and $(p+1)$\nh form bundles,
\begin{align} \label{eq:prolbnd}
  \K^p
    &= \Big( \Wedge^p \Tan^* M \oplus \Wedge^{p+1} \Tan^* M \Big)
        \otimes \V,
\end{align}
and the prolongation connection $\kcovd$ on $\K^p$, called the
\emph{Killing connection}, is given~by
\begin{align} \label{eq:prolconn}
  \kcovd_X \begin{pmatrix} \Phi \\ \Xi \end{pmatrix}
    &= \begin{pmatrix}
        \covd_X \Phi - X \intp \Xi
      \\
        \covd_X \Xi
         - \tfrac{1}{p} \Big(
            \curv_X \wedge \Phi
            - \tfrac{1}{p+1}\, X \intp (\curv \wedge \Phi)
          \Big)
      \end{pmatrix}
\end{align}
for $\Phi\in\Df^p (M, \V)$, $\Xi\in\Df^{p+1} (M, \V)$.
\begin{cor} \label{kvfprolconn}
The $\V$\nh valued Killing forms of degree $p \ge 1$ are in one-to-one
correspondence with sections of\/ $\K^p$,
\begin{align*}
  \Phi \in \Df^p (M, \V)
    &\; \leftrightarrow \;
  \Theta = \begin{pmatrix} \Phi \\ \Xi \end{pmatrix} \in \Sec(\K^p),
  & \text{where\/ $\Xi$ is given by \eqref{eq:kvfskew},}
\end{align*}
which are parallel with respect to the Killing connection $\kcovd$.
In particular, the maximal possible dimension of the solution space on a
connected manifold is $\rank \K^p = \tbinom{n+1}{p+1} \rank \V$.
\end{cor}

\subsection{Projective invariance} \label{sec:proj}
As we have already noted, the equation \eqref{eq:kvf} is projectively
invariant.
To be precise, it is invariant under a projective change of the affine
connection $\covdaff$ when considered acting on appropriately weighted
differential forms.
This is well-known in the case of the scalar-valued Killing-Yano forms.
In any case, note that the linear connection $\covdv$ on the value
bundle $\V$ must remain fixed.
The aim of this part is to compare the Killing connection in
\eqref{eq:prolconn} with the standard projective tractor connection.

Now we shall briefly recall the projective tractor calculus. For more
detail, see the references \cite{beg94}, \cite{eas08}, \cite{hsss12a},
\cite{hsss12b}, and for a more systematic approach to Cartan and
parabolic geometries we refer to the monograph \cite{cs09}, Sections
4.1.5 and 5.2.6 devoted to projective structures.
Two torsion-free affine connections $\covdaff$ and $\mcovdaff$ are
projectively equivalent if there exists a 1-form
$\Upsilon \in \Df^1(M)$ such that
\begin{align} \label{eq:projequiv}
  \mcovdaff_X Y
    & = \covdaff_X Y + \Upsilon(X) Y + \Upsilon(Y) X,
  & & \text{for all } X, Y \in \Vecf(M).
\end{align}
A \emph{projective structure} on $M$ is an equivalence class
$[\covdaff]$ of torsion-free affine connections.
The curvature of $\covdaff$ can be decomposed with respect to the
general linear group as
\begin{align} \label{eq:projcurv}
  \curvaff_{X,Y} Z
    & = \weyl_{X,Y} Z + \Rho (Y, Z) X - \Rho (X, Z) Y + \beta (X, Y) Z,
\end{align}
where $\weyl$ is the totally trace-free \emph{projective Weyl tensor},
$\Rho$ is the \emph{projective Schouten tensor} and $\beta$ is a
skew-symmetric 2-form.
Note that the projective Weyl tensor is independent of a
representative connection in the class $[\covdaff]$ and hence an invariant of the
projective structure.

We define the \emph{projective $w$-density bundles}, $w \in \er$, as
oriented line bundles
\begin{align} \label{eq:projdens}
  \Eb(w) & = \Big((\Wedge^n \Tan^* M)^{\otimes 2} \Big)^{-\frac{w}{2(n+1)}},
\end{align}
where $\Wedge^n \Tan^* M$ is the canonical line bundle of $M$.
If $\W$ is a vector bundle over $M$, we denote the
corresponding \emph{weighted bundles} $\W(w) = \W \otimes \Eb(w)$.
The affine connection $\covdaff$ canonically extends to the density
bundles as well, and a positive section $\sigma$ of $\Eb(1)$
parallel with respect to $\covdaff$ is called \emph{projective scale}.
A choice of scale $\sigma$ trivializes the density bundles, or more
generally, induces bundle isomorphisms
\begin{align} \label{eq:scaletriv}
  \W & \stackrel{\sigma}{\simeq} \W(w),
  & v & \mapsto v \otimes \sigma^w,
  & & \text{for all } v \in \Sec(W).
\end{align}
The only curvature component acting non-trivially on $\Eb(w)$
is the last one in \eqref{eq:projcurv} given by $\beta$,
\begin{align} \label{eq:projcurvdens}
  \curvaff_{X,Y} \sigma & = w \beta (X, Y) \sigma,
  & & \text{for all } \sigma \in \Sec (\Eb(w)).
\end{align}
It is convenient for our purposes to assume that $\covdaff$ is such that
$\beta = 0$, which means by \eqref{eq:projcurvdens} that $\covdaff$
admits locally a scale.
Note that this assumption is always satisfied if $\covdaff = \covdg$ is
the Levi-Civita connection of a (pseudo\nh) Riemannian metric $g$.
In this case, we have canonical global scale
\begin{align} \label{eq:gscale}
  \sigma^g & = \Big({|\det(g)|} \Big)^{-\frac{1}{2(n+1)}}
\end{align}
induced by the metric.
The existence of a scale or vanishing of $\beta$ also implies that the
Schouten tensor $\Rho$ is symmetric.

The standard \emph{projective tractor bundle} can be defined as a direct sum
\begin{align} \label{eq:stdtra}
  \Tra M & = \Tan M(-1) \oplus \Eb(-1)
\end{align}
equipped with the linear \emph{tractor connection} $\covdtra$ given by
\begin{align} \label{eq:stdtraconn}
  \covdtra_X \begin{pmatrix} Y \\ \rho \end{pmatrix}
    & = \begin{pmatrix}
          \covdaff_X Y + \rho X
        \smallskip \\
          \covdaff_X \rho - \Rho (X, Y)
        \end{pmatrix},
    & & \text{for all} \quad
        \begin{aligned}
          Y & \in \Sec (\Tan M(-1)),
        \\
          \rho & \in \Sec (\Eb(-1)).
        \end{aligned}
\end{align}
While the splitting \eqref{eq:stdtra} depends on a
representative connection in the projective class,
the tractor bundle itself and the tractor
connection are projective invariants.
The tractor connection naturally extends to other tractor bundles which
are constructed as tensors generated by $\Tra M$ and its dual
$\Tra^* M$.
In particular, the skew-symmetric tractor $(p+1)$-forms split as
\begin{align} \label{eq:traform}
  \Wedge^{p+1} \Tra^* M
    & = \big( \Wedge^p \Tan^* M \big) (p+1)
      \oplus \big( \Wedge^{p+1} \Tan^* M \big) (p+1),
\end{align}
and the tractor connection is given by
\begin{align} \label{eq:traformconn}
  \covdtra_X \begin{pmatrix} \alpha \\ \beta \end{pmatrix}
    & = \begin{pmatrix}
          \covdaff_X \alpha - X \intp \beta
        \smallskip \\
          \covdaff_X \beta + (X \intp \Rho) \wedge \alpha
        \end{pmatrix}
\end{align}
for all
$\alpha \in \Sec \big( \big( \Wedge^p \Tan^* M \big) (p+1) \big)$ and
$\beta \in \Sec \big( \big( \Wedge^{p+1} \Tan^* M \big) (p+1) \big)$.
Now we observe that in the presence of a scale we can identify the
prolongation vector bundle $\K^p$ in \eqref{eq:prolbnd} with
$\Wedge^{p+1} \Tra^*M \otimes \V$ via the isomorphisms
\eqref{eq:scaletriv}.
Moreover, we immediately see that the top slot of \eqref{eq:traformconn}
coincides with the top slot of \eqref{eq:prolconn} and corresponds to
the equation \eqref{eq:kvf} defining Killing-Yano forms in the case of
scalar-valued forms.

We proceed with the general case of vector-valued Killing forms, hence
couple the tractor connection $\covdtra$ with the connection $\covdv$ on
the vector bundle $\V$ yielding a connection $\covdtrav$ acting on
tractor-vector tensors built out of $\Tra M$ and $\V$, e.g., $\V$\nh
valued skew-symmetric tractor forms.
The formula \eqref{eq:traformconn} for the tractor connection now
becomes simply
\begin{align} \label{eq:travformconn}
  \covdtrav_X \begin{pmatrix} \Phi \\ \Xi \end{pmatrix}
    & = \begin{pmatrix}
          \covd_X \Phi - X \intp \Xi
        \smallskip \\
          \covd_X \Xi + (X \intp \Rho) \wedge \Phi
        \end{pmatrix},
\end{align}
for all
$\Phi \in \Sec \big( \big( \Wedge^p \Tan^* M \big) (p+1)
\otimes \V \big)$ and
$\Xi \in \Sec \big( \big( \Wedge^{p+1} \Tan^* M \big) (p+1)
\otimes \V \big)$.

Since the covariant derivative $\covd$ is constructed from $\covdaff$
and $\covdv$ we can split the curvature,
\begin{align}\label{eq:vfcurvsplit}
  \curv & = \curv^\covd = \curvaff + \curvv,
\end{align}
into the parts $\curvaff$ and $\curvv$ which act separately on the form part
and the value part respectively. The curvature on the form part depends
just on the affine connection $\covdaff$ and can be computed for any form
$\Phi$ as
\begin{align}\label{eq:fcurv}
  \curvaff_{X,Y} \Phi
    &= \tsum_{k=1}^n (\curvaff_{X,Y} e^k) \wedge (e_k \intp \Phi)
    = - \tsum_{k=1}^n e^k \wedge ((\curvaff_{X,Y} e_k) \intp \Phi).
\end{align}
In order to compare the curvature terms in the bottom slot of
\eqref{eq:travformconn} and \eqref{eq:prolconn} we compute the
skew-symmetrizations as defined in \eqref{eq:skewcurv1} and
\eqref{eq:skewcurv2},
\begin{align}\label{eq:skewfcurv1}
\begin{split}
  \curvaff_X \wedge \Phi
    &= - \tsum_{j,k=1}^n e^j \wedge e^k
          \wedge ((\curvaff_{X,e_j} e_k) \intp \Phi)
\\
    &= - \tfrac{1}{2} \tsum_{j,k=1}^n e^j \wedge e^k
          \wedge ((\curvaff_{X,e_j} e_k - \curvaff_{X,e_k} e_j) \intp \Phi)
\\
    &= \tfrac{1}{2} \tsum_{j,k=1}^n e^j \wedge e^k
          \wedge ((\curvaff_{e_j, e_k} X) \intp \Phi),
\end{split}
\\\label{eq:skewfcurv2}
\begin{split}
  \curvaff \wedge \Phi
    &= - \tfrac{1}{2} \tsum_{i,j,k=1}^n e^i \wedge e^j \wedge e^k
          \wedge ((\curvaff_{e_i, e_j} e_k) \intp \Phi)
\\
    &= - \tfrac{1}{6} \tsum_{i,j,k=1}^n e^i \wedge e^j \wedge e^k
          \wedge (\text{cycl.}_{(i,j,k)}(\curvaff_{e_i, e_j} e_k)
              \intp \Phi)
\\
    &= 0.
\end{split}
\end{align}
We have used repeatedly the first Bianchi identity for $\curvaff$ in the
previous computations.

\begin{rem}
In the case of scalar-valued forms we have $\curv = \curvaff$, and it is
well known that \eqref{eq:skewfcurv2} is equivalent to the first Bianchi
identity.
The completely skew-symmetrized curvature terms are the newly arising
components in the vector-valued case when compared to the case of
Killing-Yano forms, cf.\ \cite{sem03}.
\end{rem}

\begin{prop} \label{proltra}
Assume that \/ $\covdaff$ admits a scale, i.e.\ there is a \/ $\covdaff$\nh
parallel section $\sigma$ of\/ $\Eb(1)$.
Then under the identification of weighted forms with their unweighted
counterparts via the isomorphisms \eqref{eq:scaletriv}, the Killing
connection defined in \eqref{eq:prolconn} is given by
\begin{align} \label{eq:proltra}
  \kcovd_X \begin{pmatrix} \Phi \\ \Xi \end{pmatrix}
    & = \covdtrav_X \begin{pmatrix} \Phi \\ \Xi \end{pmatrix}
      - \tfrac{1}{p} \begin{pmatrix}
          0
        \\
            \weyl_X \wedge \Phi
            + \curvv_X \wedge \Phi
            - \tfrac{1}{p+1}\, X \intp (\curvv \wedge \Phi)
        \end{pmatrix},
\end{align}
for all\/ $\Phi\in\Df^p (M, \V)$ and\/ $\Xi\in\Df^{p+1} (M, \V)$.
Here the skew-symmetrized actions of\/ $\weyl$ and\/ $\curvv$ on\/
$\Phi$ are defined as in \eqref{eq:skewcurv1} and \eqref{eq:skewcurv2}.
\end{prop}

\begin{prf}
Firstly, since the scale $\sigma$ is required to be $\covdaff$\nh
parallel the isomorphisms \eqref{eq:scaletriv} preserve the covariant
derivatives and the curvature actions.
As already noted, the presence of a scale also implies $\beta = 0$ and
that the Schouten tensor $\Rho$ is symmetric.
It remains to show that the curvature terms in the bottom slot of
\eqref{eq:travformconn} and \eqref{eq:prolconn} are equal.
Indeed, the $\curvv$\nh part is left unchanged and for the $\curvaff$\nh
part we compute using \eqref{eq:skewfcurv1}, \eqref{eq:projcurv} with
$\beta = 0$, and the symmetry of $\Rho$,
\begin{align*}
\begin{split}
  \curvaff_X \wedge \Phi
    & = \weyl_X \wedge \Phi
      + \tfrac{1}{2} \tsum_{j,k=1}^n e^j \wedge e^k \wedge (
          (\Rho (e_k, X) \, e_j - \Rho (e_j, X) \, e_k) \intp \Phi
        )
\\
    & = \weyl_X \wedge \Phi
      - \tsum_{j,k=1}^n \Rho (X, e_j) \, e^j \wedge e^k \wedge (
          e_k \intp \Phi
        )
\\
    & = \weyl_X \wedge \Phi - p \, (X \intp \Rho) \wedge \Phi,
\end{split}
\end{align*}
which together with \eqref{eq:skewfcurv2} proves \eqref{eq:proltra}.
\end{prf}

Since the Weyl tensor $\weyl$ is an invariant of the projective
structure and $\covdv$ and hence also $\curvv$ are fixed independently
of\/ $\covdaff$, the equation \eqref{eq:proltra} implies that $\kcovd$
is also an invariant of the projective structure
when considered on $\V$\nh valued tractor forms.
The scale becomes redundant when appropriately weighted
forms are considered, hence by Corollary \ref{kvfprol} we get:

\begin{cor}
The Killing equation \eqref{eq:kvf} is projectively invariant when
considered on weighted $\V$\nh valued $p$-forms of weight $w = p + 1$.
\end{cor}

The corollary can be also proved directly by considering two
representative connections $\covdaff$ and $\mcovdaff$ in the projective
class and employing appropriate transformation of\/ $\Xi$ induced by the
change of splitting \eqref{eq:traform}, cf.\ also \cite{ste03}.
After some curvature rearrangements based on the first Bianchi identity,
a special case of the formula \eqref{eq:proltra} can also be found in
\cite{eas08}, p.\ 50.

\begin{rem}
A method of constructing an invariant prolongation connection via
tractors in the broad context of parabolic geometries was developed in
\cite{hsss12a}, \cite{hsss12b}.
In particular, in the Section 3.2 of \cite{hsss12b} an explicit formula
was derived for the case of skew-symmetric contravariant projective
tractors dual to our case of tractor forms.
There is a sign mistake in the relevant formula in \cite{hsss12b}, the
corrected prolongation connection $\kcovd$ is
\begin{align} \label{eq:contrtraprolidx}
  \kcovd_a
      \begin{pmatrix}
        \sigma^{\mathbf{c}}
      \\
        \rho^{\mathbf{\dot{c}}}
      \end{pmatrix}
    & = \covdtra_a
        \begin{pmatrix}
          \sigma^{\mathbf{c}}
        \\
          \rho^{\mathbf{\dot{c}}}
        \end{pmatrix}
      - \frac{\ell(\ell-1)}{2(n-\ell)}
        \begin{pmatrix}
          0
        \\
          \weyl\indices{_{pr}^{c^2}_a}
          \sigma\indices{^{pr\mathbf{\ddot{c}}}}
        \end{pmatrix}
\end{align}
and it agrees with our formula
\eqref{eq:proltra} in the special case of scalar-valued forms
($\curvv = 0$), and consequently also with \cite{sem03}.
This can be verified by a straightforward application of the bundle
isomorphism
\begin{align}
  \Wedge^\ell \Tra M \simeq \Wedge^{n+1-\ell} \Tra^* M
\end{align}
induced by a constant nonzero tractor volume form
$\Omega \in \Sec \big( \Wedge^{n+1} \Tra^* M \big)$.
We remark that it is sufficient that $\Omega$ exists at least locally,
hence the last statement holds even in the non-orientable case.

The contravariant form corresponds to the prolongation of the equation
\begin{align} \label{eq:contrkf}
  \covd_X \sigma & = - \tfrac{1}{\ell} \, X \wedge \rho,
  & & \text{for all } X \in \Vecf(M),
\end{align}
where $\sigma \in \Sec \big( \Wedge^\ell \Tan M \big)$ and
$\rho \in \Sec \big( \Wedge^{\ell-1} \Tan M \big)$ for $\ell < n$.
For $\ell = 1$ the solutions are called \emph{concircular vector
fields}, see \cite{yan40} or \cite{hbcm17}.
The equation \eqref{eq:contrkf} is clearly just a dual form of our basic
equation \eqref{eq:kvf}.
Later we introduce similar equations \eqref{eq:dkvf} and \eqref{eq:skvf}
in the presence of a (pseudo\nh) Riemannian metric.
\end{rem}

It is also worth noting that the formula \eqref{eq:prolconn} for the
Killing connection is not the only possibility to prolong the equation
\eqref{eq:kvf}.
However, Proposition \ref{proltra} implies that our particular form of
the prolongation has the advantage of being projectively invariant.

\subsection{Higher integrability conditions}

We have already proved the first integrability condition
\eqref{eq:kvfint0} in the second part of Proposition \ref{kvfprol}.
We are going to derive an explicit formula also for the second
integrability condition.
This is especially important in the case $p = 1$ when the first
condition \eqref{eq:kvfint0} becomes empty.
As it turns out the second integrability condition involves a
modification of the total curvature on $\V$\nh valued forms.

\begin{prop}\label{kvfint1}
Let\/ $\Phi$ be a Killing $\V$\nh valued form of degree $p \ge 1$ and\/
$\Xi$ the corresponding $\V$\nh valued form of degree $p+1$.
Then it holds
\begin{align}\label{eq:kvfint1}
  (\mcurv \Xi)^{(p+1,2)}_{X,Y}
    & = ((\covd \curv) \barwedge \Phi)_{X,Y},
  & & \text{for all } X,Y \in \Vecf(M),
\end{align}
where the modified curvature $\mcurv$ acting on $\Xi$ is given by
\begin{align}\label{eq:vfcurvmod}
  \mcurv &= (p+2) \curv - \curvaff
    = (p+1) \curvaff + (p+2) \curvv,
\end{align}
and the action of the derivative of the curvature on $\Phi$ is given by
\begin{align} \label{eq:barwedge}
\begin{split}
  ((\covd \curv) \barwedge \Phi)_{X,Y}
    & = (\covd_X \curv)_Y \wedge \Phi
      - (\covd_Y \curv)_X \wedge \Phi
      - {}
\\
    & \quad
      - \tfrac{1}{p+1} (Y \intp ((\covd_X \curv) \wedge \Phi)
                      - X \intp ((\covd_Y \curv) \wedge \Phi)).
\end{split}
\end{align}
The symmetry components of the action of curvature or its first
derivative are defined as in
\eqref{eq:skewcurv1}--\eqref{eq:curvdecomp3}.
\end{prop}
\begin{prf}
First we compute the action of the curvature on $\Xi$ using
\eqref{eq:kvfprol} and \eqref{eq:kvf}:
\begin{align*}
\begin{split}
  p \curv_{X,Y} \Xi
    & = \curv_Y \wedge (X \intp \Xi)
      - \curv_X \wedge (Y \intp \Xi)
      - {}
\\
  & \quad
      - \tfrac{1}{p+1} (Y \intp (\curv \wedge (X \intp \Xi))
                      - X \intp (\curv \wedge (Y \intp \Xi)))
      + {}
\\
  & \quad
      + (\covd_X \curv)_Y \wedge \Phi
      - (\covd_Y \curv)_X \wedge \Phi
      - {}
\\
  & \quad
      - \tfrac{1}{p+1} (Y \intp ((\covd_X \curv) \wedge \Phi)
                      - X \intp ((\covd_Y \curv) \wedge \Phi))
\end{split}
\end{align*}
In the next we compute the terms containing the contraction of $\Xi$
using \eqref{eq:skewcurv1} and~\eqref{eq:skewcurv2}:
\begin{align*}
\begin{split}
  \curv_Y \wedge (X \intp \Xi)
    & = \tsum_{j=1}^n e^j \wedge ((\curvaff_{Y,e_j} X) \intp \Xi
                              + X \intp (\curv_{Y,e_j} \Xi))
\\
  & = \tsum_{j=1}^n e^j \wedge ((\curvaff_{Y,e_j} X) \intp \Xi)
      - \curv_{X,Y} \Xi - X \intp (\curv_Y \wedge \Xi),
\end{split}
\\
\begin{split}
  \curv \wedge (X \intp \Xi)
    & = \tfrac{1}{2}\tsum_{i,j=1}^n e^i \wedge e^j
                            \wedge ((\curvaff_{e_i,e_j} X) \intp \Xi)
      - \curv_X \wedge \Xi + X \intp (\curv \wedge \Xi) .
\end{split}
\end{align*}
Collecting all the terms containing $\Xi$ and using \eqref{eq:fcurv},
\eqref{eq:skewfcurv1}, \eqref{eq:skewfcurv2} and the first Bianchi
identity (recall that $\curvaff$ acts just on the form part), we get
\begin{align*}
\begin{split}
  & p \curv_{X,Y} \Xi
      - \curv_Y \wedge (X \intp \Xi)
      + \curv_X \wedge (Y \intp \Xi)
      + {}
\\
  &
      + \tfrac{1}{p+1} (Y \intp (\curv \wedge (X \intp \Xi))
                      - X \intp (\curv \wedge (Y \intp \Xi)))
    = {}
\end{split}
\\
\begin{split}
  & \quad
    = (p+2) \Big(
        \curv_{X,Y} \Xi
        - \tfrac{1}{p+1} \Big( Y \intp (\curv_X \wedge \Xi)
                            - X \intp (\curv_Y \wedge \Xi) - {}
\\
  & \qquad \qquad \qquad \qquad \qquad \qquad \;
          - \tfrac{2}{p+2}\, Y \intp (X \intp (\curv \wedge \Xi))
        \Big)
      \Big)
      - {}
\\
  & \quad \quad
      - \Big(
        \curvaff_{X,Y} \Xi
        - \tfrac{1}{p+1} \Big( Y \intp (\curvaff_X \wedge \Xi)
                            - X \intp (\curvaff_Y \wedge \Xi) - {}
\\
  & \qquad \qquad \qquad \qquad \qquad \,
          - \tfrac{2}{p+2}\, Y \intp (X \intp (\curvaff \wedge \Xi))
        \Big)
      \Big).
\end{split}
\end{align*}
Finally \eqref{eq:kvfint1} follows by \eqref{eq:curvdecomp3} for
$\Xi$ (recall that the degree of $\Xi$ is $p+1$).
\end{prf}

In general, there are many other integrability conditions for
Killing $\V$\nh valued forms resulting from the prolongation
procedure as formulated in Corollary \ref{kvfprolconn}.
The complete set of integrability conditions for the equation
$\kcovd \Theta = 0$ is given as the annihilator of the infinitesimal
holonomy algebra $\hola(\K^p, \kcovd)$.
The infinitesimal holonomy algebra is pointwise generated by curvature
and all its derivatives, hence the integrability conditions are
\begin{align}\label{eq:parint}
  (\kcovd^k_{Z_1,\dots,Z_k} \kcurv)_{X,Y} \Theta &= 0,
  & \begin{array}{l}
      \text{for all } X,Y,Z_1,\dots,Z_k \in \Vecf(M)\,
    \\
      \text{and } k = 0, 1, \dots .
    \end{array}
\end{align}
In other words, the integrability conditions reflect the obstructions to
flatness of the prolongation tractor bundle with respect to the Killing
connection.

However, we find it more convenient for further application to formulate
the integrability conditions directly in terms of the components $\Phi$,
$\Xi$ as in the formulas \eqref{eq:kvfint0} and \eqref{eq:kvfint1}.
This corresponds to splitting \eqref{eq:parint} into individual slots.
In fact, equation \eqref{eq:kvfint0} appears in the top slot of
\eqref{eq:parint} for $k = 0$, and equation \eqref{eq:kvfint1} appears
in the bottom slot of \eqref{eq:parint} for $k = 0$ as well as in an
equivalent form in the top slot for $k = 1$.

\section{Special Killing vector-valued differential forms}
\label{sec:skvf}

From now on $(M, g)$ is assumed to be a (pseudo\nh) Riemannian manifold
and the affine connection $\covdaff$ will always be the
Levi-Civita connection $\covdaff = \covdg$.
We will denote the corresponding (pseudo\nh) Riemannian curvature by
$\curvg$.

The curvature operator at a point $x\in M$ takes values in the
orthogonal Lie algebra $\sola(\Tan_x M, g)$, identified with
$\Wedge^2 \Tan_x M$ by an isomorphism $\rho$:
\begin{align} \label{eq:sola}
  \rho(X \wedge Y) \, Z &= g(Y,Z) \, X - g(X,Z) \, Y,
  & & \text{for all } X, Y, Z \in \Tan_x M.
\end{align}
The action on skew-symmetric forms is induced by the tensor product
action,
\begin{align} \label{eq:solaform}
\begin{split}
  \rho(X \wedge Y) \, \Phi
    & = X \intp (Y^\flat \wedge \Phi) - Y \intp (X^\flat \wedge \Phi)
\\
    & = X^\flat \wedge (Y \intp \Phi) - Y^\flat \wedge (X \intp \Phi),
\end{split}
  & \begin{array}{l}
      \text{for all } X, Y \in \Tan_x M,
    \\
      \Phi \in \Wedge^p \Tan^*_x M,
    \end{array}
\end{align}
where $X^\flat$ denotes the metric dual of the vector $X$ with respect
to $g$.
The presence of metric allows us to dualize the notion of a Killing
form.

\begin{dfn} \label{dkvf}
Let $\Xi$ be a $\V$\nh valued form of degree $p \in \{1, \dots, n\}$.
We say that $\Xi$ is a \emph{$\star$\nh Killing $\V$\nh valued form} if
there exists a $\V$\nh valued form $\Phi$ of degree $p-1$ such that
\begin{align} \label{eq:dkvf}
  \covd_X \Xi &= X^\flat \wedge \Phi,
  & & \text{for all } X \in \Vecf(M).
\end{align}
\end{dfn}

The $\star$\nh Killing forms are just the Hodge star duals of Killing forms.
An interested reader can find more detailed treatment of the scalar-valued case
in \cite{sem03}.
There is a special important class of Killing forms given by matching pairs
of a Killing and a $\star$\nh Killing form. Note that in
the following definition we allow the (non-trivial) case $p = 0$.

\begin{dfn} \label{skvf}
Let $\Phi$ be a Killing $\V$\nh valued form of degree
$p \in \{0, \dots, n\}$ and $\Xi$ the corresponding $\V$\nh valued form
of degree $p+1$.
We say that $\Phi$ is a \emph{special Killing $\V$\nh valued form} if
there exists a constant $c \in \er$ such that in addition to the
equation \eqref{eq:kvf} it holds
\begin{align} \label{eq:skvf}
  \covd_X \Xi &= - c X^\flat \wedge \Phi,
  & & \text{for all } X \in \Vecf(M).
\end{align}
\end{dfn}

The main significance of special Killing forms stems from their
close relationship with the \emph{metric cone construction}, see
Section \ref{sec:cone} for more details. We mention an example
of \emph{Sasakian structures}, equivalent to special Killing-Yano
$1$-forms of constant length $1$, cf.\ \cite{sem03}.

We note that the system of equations \eqref{eq:kvf} and \eqref{eq:skvf}
is already in a closed form and there is no need to prolong it. In the
following proposition we give the first integrability condition which in
fact characterizes special Killing forms among Killing forms.

\begin{prop} \label{skvfint0}
Let\/ $\Phi$ be a Killing $\V$\nh valued form.
Then it is a special Killing $\V$\nh valued form with the corresponding
constant $c \in \er$ if and only if it holds
\begin{align} \label{eq:skvfint0}
  \curv_{X,Y} \Phi & = c \rho^\aff(X \wedge Y) \, \Phi,
  & & \text{for all } X, Y \in \Vecf(M),
\end{align}
where the superscript `$\aff$' emphasizes that it acts only on the form
part of\/ $\Phi$.
\end{prop}

\begin{prf}
First suppose that $\Phi$ is a special Killing $\V$\nh valued form of
degree $p$ and $\Xi$ the corresponding $\V$\nh valued form of degree
$p+1$.  Using \eqref{eq:kvf}, \eqref{eq:skvf} and \eqref{eq:solaform} we
compute the action of the curvature proving \eqref{eq:skvfint0}:
\begin{align*}
\begin{split}
  \curv_{X,Y} \Phi
    & = \covd^2_{X,Y} \Phi - \covd^2_{Y,X} \Phi
    = Y \intp (\covd_X \Xi) - X \intp (\covd_Y \Xi)
\\
    & = -c (Y \intp (X^\flat \wedge \Phi)
          - X \intp (Y^\flat \wedge \Phi))
    = c \rho^\aff(X \wedge Y)\, \Phi.
\end{split}
\end{align*}

On the other hand suppose that $\Phi$ is Killing and \eqref{eq:skvfint0}
holds.
We will compute in a tangent frame $\{e_1,\ldots ,e_n\}$, with
$g_{ij} = g(e_i, e_j)$ the corresponding metric components.
The skew-symmetrized actions of the curvature as defined in
\eqref{eq:skewcurv1} and \eqref{eq:skewcurv2} are
\begin{align*}
\begin{split}
  \curv_X \wedge \Phi
    & = c \tsum_{j=1}^n e^j \wedge (\rho^\aff(X \wedge e_j) \, \Phi)
\\
    & = c \tsum_{j=1}^n e^j \wedge (X^\flat \wedge (e_j \intp \Phi)
                                  - (e_j)^\flat \wedge (X \intp \Phi))
\\
    & = - c \Big( X^\flat
                \wedge \Big(
                    \tsum_{j=1}^n e^j \wedge (e_j \intp \Phi)
                  \Big)
          + \tsum_{i,j=1}^n g_{ij} \, e^j \wedge e^i
                                      \wedge (X \intp \Phi)
        \Big)
\\
    & = - c p X^\flat \wedge \Phi,
\end{split}
\\
  \curv \wedge \Phi
    & = - \tfrac{cp}{2} \tsum_{i=1}^n e^i \wedge (e_i)^\flat \wedge \Phi
    = - \tfrac{cp}{2} \tsum_{i,j=1}^n g_{ij} \, e^i \wedge e^j
                                                    \wedge \Phi
    = 0.
\end{align*}
Now \eqref{eq:skvf}
follows by using Proposition \ref{kvfprol} and substituting into
\eqref{eq:kvfprol}.
\end{prf}

We conclude this section with a few higher integrability conditions.

\begin{prop} \label{skvfint1}
Let $\Phi$ be a special Killing $\V$\nh valued $p$\nh form, $\Xi$ the
corresponding $\V$\nh valued $(p+1)$\nh form and $c \in \er$ the
associated Killing constant.
Then the following equalities hold for all $X, Y, Z \in \Vecf(M)$:
\begin{align} \label{eq:skvfint1}
  \curv_{X,Y} \Xi & = c \rho^\aff(X \wedge Y) \, \Xi,
\medskip
\\ \label{eq:skvfint2}
  (\covd_X \curv)_{Y,Z} \Phi
    & = -((\curvg_{Y,Z} - c \rho(Y \wedge Z)) \, X) \intp \Xi,
\medskip
\\ \label{eq:skvfint3}
  (\covd_X \curv)_{Y,Z} \Xi
    & = c ((\curvg_{Y,Z} - c \rho(Y \wedge Z)) \, X)^\flat \wedge \Phi.
\end{align}
\end{prop}

\begin{prf}
The proof of \eqref{eq:skvfint1} is analogous to that of
\eqref{eq:skvfint0}, we just compute the action of the curvature using
\eqref{eq:kvf}, \eqref{eq:skvf} and \eqref{eq:solaform}:
\begin{align*}
\begin{split}
  \curv_{X,Y} \Xi
    & = \covd^2_{X,Y} \Xi - \covd^2_{Y,X} \Xi
    = -c (Y^\flat \wedge (\covd_X \Phi)
        - X^\flat \wedge (\covd_Y \Phi))
\\
    & = -c (Y^\flat \wedge (X \intp \Xi)
        - X^\flat \wedge (Y \intp \Xi))
    = c \rho^\aff(X \wedge Y)\, \Xi.
\end{split}
\end{align*}
Equations \eqref{eq:skvfint2} and \eqref{eq:skvfint3} follow from
\eqref{eq:skvfint0} and \eqref{eq:skvfint1} by differentiating while
noticing that $\rho(\cdot \wedge \cdot)$ is covariantly constant.
\end{prf}

\begin{rem}
The integrability conditions can be elegantly formulated in terms of a
curvature modified by the scalar curvature component in the Ricci
decomposition,
\begin{align} \label{eq:curvconc}
  \curv^c_{X,Y} & := \curv_{X,Y} - c \rho^\aff(X \wedge Y).
\end{align}
The second term is covariantly constant and thus we have
$\covd \curv = \covd \curv^c$, so the modified curvature can be
conveniently incorporated into all the higher integrability
conditions.

This is related to the \emph{concircular curvature} consisting just of
the trace-free Ricci and (conformal) Weyl components, cf.\ \cite{yan40}.
The modification has also nice interpretation from the point of view
of Cartan geometries and the phenomenon of \emph{model mutation}, cf.\
\cite{cs09} Sections 1.1.2 and 1.5.1.
It is the modification we would get if we had chosen the round sphere
(for $c > 0$) or the hyperbolic space (for $c < 0$), respectively, as
the homogeneous model space for the Riemannian geometry instead of the
Euclidean space (and analogously in the pseudo-Riemannian~case).
\end{rem}

\section{Killing spinor-valued differential forms}
\label{sec:ksf}

From now on $(M, g)$ is assumed to be oriented and spin with a fixed
chosen spin structure.
We denote the corresponding complex spinor bundle by $\Spnr M$ and the
operation of Clifford multiplication by vectors on spinors by `$\clp$'.
The Levi-Civita connection $\covdg$ lifts to a spin connection on
$\Spnr M$, and we denote it by abuse of notation $\covdg$ too.
We briefly recall the notion of Killing spinors and their basic
properties relevant for our needs, see for example
the monograph \cite{bfgk91}.

\begin{dfn} \label{ks}
A spinor field $\Psi$ is a \emph{Killing spinor} if there exists a
constant $a \in \ce$ such that
\begin{align} \label{eq:ks}
  \covdg_X \Psi &= a X \clp \Psi,
  & & \text{for all } X \in \Vecf(M).
\end{align}
The constant $a$ is called the \emph{Killing number} of $\Psi$.
\end{dfn}

The equation \eqref{eq:ks} is already in the closed form and we can
readily write down the prolongation connection, called again the
\emph{Killing connection} on spinors,
\begin{align}\label{eq:covda}
  \covda_X \Psi & = \covdg_X \Psi - a X \clp \Psi,
  & X \in \Vecf(M), \, \Psi \in \Sec (\Spnr M).
\end{align}

\begin{cor} \label{ksprol}
The Killing spinors with fixed Killing number $a \in \ce$ are
sections of\/ $\Spnr M$ parallel with respect to the Killing
connection $\covda$.
In particular, the maximal possible dimension of the solution space on a
connected manifold is $\rank \Spnr M = 2^{\lfloor \tfrac{n}{2} \rfloor}$.
\end{cor}

By this corollary or directly from \eqref{eq:ks} it follows easily the
first integrability condition.
We recall the action of the Lie algebra
$\spinla (\Tan_x M, g) \simeq \sola (\Tan_x M, g)$ on spinors in terms of
the isomorphism $\rho$ in \eqref{eq:sola},
\begin{align} \label{eq:spinla}
  \rho(X \wedge Y) \, \Psi
    &= - \tfrac{1}{4} [{X \clp}, {Y \clp}] \, \Psi
    = - \tfrac{1}{2} (X \clp {Y \clp} + g(X, Y)) \, \Psi,
\end{align}
where $[\cdot, \cdot]$ denotes the commutator.
The curvature of the Killing connection $\covda$ defined in
\eqref{eq:covda} is given by
\begin{align} \label{eq:curva}
\begin{split}
  \curva_{X,Y} \Psi
    & = \covda_X (\covda_Y \Psi) - \covda_Y (\covda_X \Psi)
      - \covda_{[X,Y]} \Psi
\\
    & = \curvg_{X,Y} \Psi + a^2 [{X \clp}, {Y \clp}] \, \Psi
    = \curvg_{X,Y} \Psi - 4a^2 \rho(X \wedge Y) \, \Psi.
\end{split}
\end{align}

\begin{prop} \label{ksint0}
Let $\Psi$ be a Killing spinor with Killing number $a \in \ce$. Then
\begin{align} \label{eq:ksint0}
  \curvg_{X,Y} \Psi
    & = -a^2 \, [{X \clp}, {Y \clp}] \, \Psi
    = 4 a^2 \rho(X \wedge Y) \, \Psi,
  & & \text{for all } X, Y \in \Vecf(M).
\end{align}
\end{prop}

By \eqref{eq:ksint0} it follows the well-known fact that a manifold
admitting Killing spinors has constant \emph{scalar curvature},
related to the Killing number by
\begin{align} \label{eq:ksscal}
  \Scal^g & = 4 n (n - 1) a^2.
\end{align}
In the Riemannian case it also follows that the manifold is Einstein.

We define the Killing spinor-valued forms as a special case of Killing
vector-valued forms, where $\covdaff = \covdg$ as before and the vector
bundle of values is the spinor bundle $\V = \Spnr M$ equipped with the
connection $\covdv = \covda$ given above for arbitrary $a \in \ce$.
We can now reformulate Definitions \ref{kvf} and \ref{skvf} in terms
of the Levi-Civita spin connection $\covdg$ and the Killing number
$a \in \ce$.
In particular, by abuse of notation $\covdg$ denotes the usual tensor
product connection acting on both the form and the spinor part of a
spinor-valued form by $\covdg$.

\begin{dfn} \label{ksf}
Let $\Phi$ be a spinor-valued form of degree $p \in \{0, \dots, n\}$.
We say that $\Phi$ is a \emph{Killing spinor-valued form} if there
exists a constant $a \in \ce$ and a spinor-valued form $\Xi$ of degree
$p+1$ such that
\begin{align} \label{eq:ksf}
  \covdg_X \Phi & = a X \clp \Phi + X \intp \Xi,
  & & \text{for all } X \in \Vecf(M).
\end{align}
We say that $\Phi$ is a \emph{special Killing spinor-valued form} if
there exists another constant $c \in \er$ such that both \eqref{eq:ksf} and
the equation
\begin{align} \label{eq:sksf}
  \covdg_X \Xi & = a X \clp \Xi - c X^\flat \wedge \Phi,
  & & \text{for all } X \in \Vecf(M)
\end{align}
are satisfied.
\end{dfn}

In order to develop calculus for spinor-valued forms it is convenient to express
the Clifford multiplication ${\clp} \colon \Tan M \otimes \Spnr M \to \Spnr
M$ as a 1-form with values in the endomorphisms of the spinor bundle
$\clf \in \Omega^1(M, \End(\Spnr M))$,
\begin{align} \label{eq:clf}
  \clf &= \tsum_{i=1}^n e^i \otimes ({e_i \clp}).
\end{align}
In terms of $\clf$ the defining equation of the Clifford algebra can be
expressed as
\begin{align} \label{eq:cla}
  \text{symm.}(\clf \otimes \clf) & = -2 g,
\end{align}
where $\text{symm.}$ denotes the symmetrization in form indices. Note that the
Clifford-valued form $\clf$ is parallel with respect to the Levi-Civita
spin connection $\covdg$ because the Clifford multiplication is
equivariant with respect to the spin group.  We will also frequently use
the following algebraic formulas,
\begin{align} \label{eq:clf1}
  X \clp (\clf \wedge \Phi) + \clf \wedge (X \clp \Phi)
    & = -2 X^\flat \wedge \Phi,
\\ \label{eq:clf2}
  X \intp (\clf \wedge \Phi) + \clf \wedge (X \intp \Phi)
    & = X \clp \Phi,
\\ \label{eq:clf3}
  X \clp (\clv \intp \Phi) + \clv \intp (X \clp \Phi)
    & = -2 X \intp \Phi,
\\ \label{eq:clf4}
  X^\flat \wedge (\clv \intp \Phi) + \clv \intp (X^\flat \wedge \Phi)
    & = X \clp \Phi,
\end{align}
as a consequence of \eqref{eq:clf} and \eqref{eq:cla}.
Here
\begin{align} \label{eq:clv}
  \clv & = \tsum_{i,j=1}^n g^{ij} \, e_i \otimes ({e_j \clp})
\end{align}
is the metric dual of the 1-form $\clf$ with respect to $g$ and
$g^{ij} = g(e^i, e^j) = (g_{ij})^{-1}$ are the entries of the inverse
metric.
By \eqref{eq:clf}--\eqref{eq:clv},
\begin{align} \label{eq:clfclv}
  \clv \intp (\clf \wedge \Phi) - \clf \wedge (\clv \intp \Phi)
    & = (2 p - n) \, \Phi,
\end{align}
for all spinor-valued forms $\Phi$ of degree $p$.

In the present case the equation \eqref{eq:kvfskew} which determines
$\Xi$ specializes to
\begin{align} \label{eq:ksfskew}
  \Xi & = \tfrac{1}{p+1} \, \difa \Phi
    = \tfrac{1}{p+1} (\difg \Phi - a \clf \wedge \Phi),
\end{align}
where $\difa$ and $\difg$ denote the respective exterior covariant
derivatives with $\covda$ and $\covdg$ acting on the spinor values.
We can now substitute $\Xi$ into \eqref{eq:ksf} to get another
equivalent form of the defining equation,
\begin{align} \label{eq:ksfsub}
  \covdg_X \Phi
    & = a \big( X \clp \Phi
        - \tfrac{1}{p+1} \, X \intp (\clf \wedge \Phi)
      \big)
      + \tfrac{1}{p+1} \, X \intp \difg \Phi.
\end{align}
Similarly we can substitute $\Xi$ into \eqref{eq:sksf}, and by \eqref{eq:clf1}, \eqref{eq:clf2} get
\begin{align} \label{eq:sksfsub}
\begin{split}
  & \covdg_X (\difg \Phi)
    = \covdg_X ((p+1) \, \Xi + a \clf \wedge \Phi)
\\ & \qquad
    = (p+1) (a X \clp \Xi - c X^\flat \wedge \Phi)
      + a \clf \wedge (a X \clp \Phi + X \intp \Xi)
\\ & \qquad
    = -c (p+1) \, X^\flat \wedge \Phi
      + a \Big( X \clp \difg \Phi
              + \tfrac{1}{p+1} \, \clf \wedge (X \intp \difg \Phi)
            \Big)
      + {}
\\ & \qquad \quad
      + a^2 \Big( 2 X^\flat \wedge \Phi
              + \tfrac{2p+1}{p+1} \, \clf \wedge (X \clp \Phi)
              + \tfrac{1}{p+1} \, \clf \wedge (\clf
                                        \wedge (X \intp \Phi))
            \Big).
\end{split}
\end{align}
The resulting formulas \eqref{eq:ksfsub} and \eqref{eq:sksfsub} agree
with the definitions of (special) Killing spinor-valued forms introduced
in \cite{sz16}.

\section{Cone construction}
\label{sec:cone}

We shall briefly recall the cone construction following the reference
\cite{sz16}, which
provides a useful description of special Killing spinor-valued forms.
As it turns out in Section \ref{sec:ccurv}, most of the
Killing spinor-valued forms on spaces of constant curvature are special
and hence arise from this construction.

Let $(M, g)$ be a spin (pseudo\nh) Riemannian manifold of signature
$(n_+, n_-)$ with a fixed spin structure, and for $\eps = \pm 1$ we define
\begin{align} \label{eq:conesig}
\begin{aligned}
  \cone{n}_+ & = n_+ + 1, & \qquad \cone{n}_- & = n_-,
  & \qquad & \text{if } \eps = +1,
\\
  \cone{n}_+ & = n_+, &\qquad  \cone{n}_- & = n_- + 1,
  & \qquad & \text{if } \eps = -1.
\end{aligned}
\end{align}
The \emph{$\eps$-metric cone} over $(M,g)$ is the product manifold
$\cone{M} = \er_+ \times M$ with the warped product metric $\cone{g}$ of
signature $(\cone{n}_+, \cone{n}_-)$,
\begin{align} \label{eq:coneg}
  \cone{g} = \eps \dr^2 + r^2 g,
\end{align}
where $r$ is the coordinate function on $\er_+$.
The original manifold $(M, g)$ is isometrically embedded in
$(\cone{M}, \cone{g})$ as the hypersurface defined by $r = 1$ and the
outer unit normal is given by the radial vector field $\vr$.
Note that $\cone{g}(\vr, \vr) = \eps$.

\subsection{Tangent bundle and forms on hypersurfaces}

Let us consider the restricted vector bundle
$\cone{\Tan} M = \Tan \cone{M} |_M$ over $M$.
It splits orthogonally into the normal and tangent bundle of $M$,
\begin{align}
  \cone{\Tan} M & = \Norm M \oplus \Tan M.
\end{align}
The normal bundle $\Norm M$ is in our case trivialized by the outer unit
normal $\vr$.
Hence we can write the normal and tangent projections as
\begin{align}
\begin{split}
  \pi^\Norm (\cone{v}) & = \cone{g} (\vr, \cone{v}),
\\
  \pi^\Tan (\cone{v})
    & = \cone{v} - \eps \cone{g} (\vr, \cone{v}) \, \vr,
\end{split}
  & & \text{for all } \cone{v} \in \cone{\Tan} M.
\end{align}
Accordingly, the decomposition of differential forms is
\begin{align} \label{eq:conefdecomp}
  \Wedge^{p+1} \cone{\Tan}^* M
    & = (\Norm^* M \otimes \Wedge^p \Tan^* M)
      \oplus \Wedge^{p+1} \Tan^* M,
\end{align}
and the corresponding projections are given by
\begin{align} \label{eq:conefproj}
\begin{split}
  \pi^\Norm (\cone{\alpha})
    & = \vr \intp \cone{\alpha},
\\
  \pi^\Tan (\cone{\alpha})
    & = \cone{\alpha} - \dr \wedge (\vr \intp \cone{\alpha}),
\end{split}
  & & \text{for all } \cone{\alpha} \in \Wedge^{p+1} \cone{\Tan}^* M.
\end{align}

The \emph{shape operator} of $M$ regarded as embedded in its $\eps$-metric cone is
$\Shape(X) = -X$ and hence the respective Levi-Civita connections
$\covdg$ and $\cone{\covdg}$ are related by
\begin{align} \label{eq:conecovdg}
  \covdg_X Y & = \cone{\covdg_X} Y + \eps g(X, Y) \, \vr,
  & & \text{for all } X, Y \in \Vecf(M),
\end{align}
cf.\ the formulas for $\cone{\covdg}$ in \cite{sem03}, \cite{sz16}.
It is convenient to extend $\covdg$ to the whole vector bundle
$\cone{\Tan} M$ so that the normal vector field $\vr$ is covariantly constant, which
yields the general formula
\begin{align} \label{eq:econecovdg}
  \covdg_X
    & = \cone{\covdg_X} + \eps \cone{\rho}(\vr \wedge X),
  & & \text{for all } X \in \Vecf(M),
\end{align}
where
$\cone{\rho} \colon \Wedge^2 \cone{\Tan}_x M
  \to \sola(\cone{\Tan}_x M, \cone{g})$
is defined analogously as in \eqref{eq:sola}.
In particular, using \eqref{eq:solaform} we get for skew-symmetric forms
\begin{align} \label{eq:econecovdgf}
  \covdg_X \cone{\alpha}
    & = \cone{\covdg_X} \cone{\alpha}
      + \dr \wedge (X \intp \cone{\alpha})
      - \eps \, X^\flat \wedge (\vr \intp \cone{\alpha}),
\end{align}
for all $X \in \Vecf(M)$,
$\cone{\alpha} \in \Sec(\Wedge^p \cone{\Tan}^* M)$.
By construction, the extended connection $\covdg$ commutes with the
projections $\pi^\Norm$ and $\pi^\Tan$.

\subsection{Spinors on hypersurfaces}

The description of spinor bundles on a hypersurface can be found in
\cite{bar96}, \cite{bgm05}.
Here we discuss the case of our interest for signature $(n_+,n_-)$ and
space- or time-like normal depending on $\eps = \pm 1$.

Firstly, we can naturally realize the Clifford algebra $\Cl(n, \ce)$ as a
subalgebra of $\Cl(n+1, \ce)$, and the corresponding complex spinor
space $\Spnr_n$ as a $\Cl(n, \ce)$-submodule of $\Spnr_{n+1}$.
We also recall that $\Cl(n, \ce)$ is isomorphic to the even
part $\Cl_0(n+1, \ce)$ of $\Cl(n+1, \ce)$.
In particular, there are two such isomorphisms
$\varphi_\pm \colon \Cl(n, \ce) \isoto \Cl_0(n+1, \ce)$
given on the generators by
\begin{align} \label{cl+1}
  \varphi_\pm (e_i) & = \mp \sqrt{\eps} \, e_0 \clp e_i,
  & & i = 1, \dots, n,
\end{align}
where $\{e_0, \dots, e_n\}$ is an orthonormal basis of
$(\er^{n+1}, \cone{g})$ such that $\cone{g}(e_0, e_0) = \eps$.
All formulas are valid for both choices of the square root sign and we
fix $\sqrt{\eps} = 1$ for $\eps = 1$ and $\sqrt{\eps} = \ii$ for
$\eps = -1$, respectively.
There are also corresponding endomorphisms
$f_\pm \colon \Spnr_{n+1} \to \Spnr_{n+1}$,
\begin{align} \label{spnr+1}
  f_\pm (\cone{\psi}) & = (1 \mp \sqrt{\eps} \, e_0) \clp \cone{\psi},
  & & \text{for all } \cone{\psi} \in \Spnr_{n+1},
\end{align}
which intertwine the restricted representations of $\Cl(n, \ce)$ and
$\Cl_0(n+1, \ce)$ on $\Spnr_{n+1}$ with respect to the isomorphisms
$\varphi_\pm$.  The mappings $f_+$ and $f_-$ are up to a scalar multiple
mutual inverses,
\begin{align} \label{eq:spnr+1inv}
  f_\pm \circ f_\mp & = (1 - \eps \, e_0^2) = 2,
\end{align}
and so are linear isomorphisms.
Restricting $f_\pm$ to the subspace $\Spnr_n$ we get
\begin{align} \label{eq:spnr+1spc}
\begin{aligned}
  \Spnr_{n+1} & = f_+ (\Spnr_n) = f_- (\Spnr_n),
  & & \text{ if $n$ is even},
\\
  \Spnr_{n+1} & = f_+ (\Spnr_n) \oplus f_- (\Spnr_n),
  & & \text{ if $n$ is odd},
\end{aligned}
\end{align}
as $\Cl_0(n+1, \ce)$-modules.
In particular, the odd case agrees with the decomposition of
$\Spnr_{n+1}$ to \emph{half-spinors}, so adopting suitable sign
convention we have
\begin{align} \label{eq:hspnr+1}
  \Spnr^+_{n+1} & = f_+ (\Spnr_n),
  & \Spnr^-_{n+1} & = f_- (\Spnr_n).
\end{align}

Recall that we assume a fixed spin structure on $M$ and since the
$\eps$-metric cone $\cone{M} = \er_+ \times M$ is homotopy equivalent to
$M$, there is a unique compatible spin structure on $\cone{M}$.
Now we pass to the associated complex spinor bundles $\Spnr M$ and
$\Spnr \cone{M}$ and denote the restricted bundle
$\cone{\Spnr} M = \Spnr \cone{M} |_M$.
Compatibility of the spin structures implies that $\Spnr M$ is naturally
a subbundle of $\cone{\Spnr} M$.
The linear automorphisms $f_\pm$ induce bundle automorphisms
$F_\pm \colon \cone{\Spnr} M \to \cone{\Spnr} M$,
\begin{align} \label{eq:F}
  F_\pm (\cone{\Psi}) & = (1 \mp \sqrt{\eps} \, \vr) \clp \cone{\Psi},
  & & \text{for all } \cone{\Psi} \in \cone{\Spnr} M.
\end{align}
Since $f_\pm$ are intertwining with respect to $\varphi_\pm$, we have
the identity
\begin{align} \label{eq:Fclp}
  F_\pm (X \clp \cone{\Psi})
    & = \mp \sqrt{\eps} \, \vr \clp X \clp F_\pm (\cone{\Psi}),
  & & \text{for all } X \in \Vecf(M), \,
        \cone{\Psi} \in \Sec(\cone{\Spnr} M).
\end{align}
The formula \eqref{eq:econecovdg} together with \eqref{eq:spinla} yields
the spin connection $\covdg$ on $\cone{\Spnr} M$,
\begin{align} \label{eq:econecovdspin}
  \covdg_X \cone{\Psi}
    & = \cone{\covdg_X} \cone{\Psi}
        - \tfrac{1}{2} \, \eps \, \vr \clp X \clp \cone{\Psi},
    & & \text{for all } X \in \Vecf(M), \,
        \cone{\Psi} \in \Sec(\cone{\Spnr} M).
\end{align}
Finally, we can obtain the spin connection $\covdg$ on $\Spnr M$ by
pulling back along the bundle map $F_+$ or $F_-$ and restricting to
$\Spnr M$.
Due to $\covdg \vr = 0$, we have $\covdg F_\pm = 0$ and so
$\covdg$ commutes with $F_\pm$, hence the formula \eqref{eq:econecovdspin}
remains unchanged after the pull-back.

\subsection{Equivalences for Killing equations}

The decomposition \eqref{eq:conefdecomp} suggests that $(p+1)$-forms
over the extended tangent bundle $\cone{\Tan} M$ are isomorphic to the
prolongation vector bundle $\K^p$ defined in~\eqref{eq:prolbnd}.
Moreover, if we decompose the Levi-Civita connection $\cone{\covdg}$ on
the cone into individual slots, cf.\ \eqref{eq:econecovdgf}, we basically
obtain the defining equations \eqref{eq:kvf} and \eqref{eq:skvf} for
\emph{special} Killing-Yano forms.
Hence we arrive at the following equivalence, see \cite{sem03} or
\cite{sz16} for a detailed proof.

\begin{prop} \label{coneskf}
Let\/ $\alpha$ and\/ $\beta$ be differential forms on $M$ of degrees $p$
and $p+1$ respectively, and define a differential form\/ $\cone{\theta}$
on the $\eps$-metric cone $\cone{M}$ by
\begin{align} \label{eq:coneskf}
  \cone{\theta}
    & = r^p \dr \wedge \pi_2^* (\alpha) + r^{p+1} \pi_2^* (\beta),
\end{align}
where $\pi_2^*$ denotes the pull-back along the canonical projection
$\pi_2 \colon \cone{M} \to M$.
Then\/ $\cone{\theta}$ is parallel with respect to\/ $\cone{\covdg}$ if
and only if\/ $\alpha$ is special Killing-Yano form with the
corresponding $(p+1)$-form\/ $\beta$ and the Killing constant
$c = \eps$.

Conversely, any parallel differential form\/ $\cone{\theta}$ of degree
$p+1$ on $\cone{M}$ arises this way with\/ $\alpha$ and\/ $\beta$ given
by the normal and tangent projections respectively,
\begin{align} \label{eq:coneskfinv}
  \alpha & = \pi^\Norm (\cone{\theta} |_M),
  & \beta & = \pi^\Tan (\cone{\theta} |_M).
\end{align}
\end{prop}

The homogeneity factors $r^p$ and $r^{p+1}$ in \eqref{eq:coneskf} ensure
that $\cone{\theta}$ is parallel in the direction $\vr$.
This is equivalent to the projective weight $w = p + 1$ which appears in
\eqref{eq:traform}.
There is a related construction of the so called \emph{Thomas projective
cone} equipped with an affine connection which is equivalent to the
standard projective tractor connection, see \cite{arm08}.
In particular, for Einstein manifolds the two cone constructions
essentially coincide, and special Killing-Yano forms are thus equivalent
to parallel tractor forms.
This approach is further exploited~in~\cite{gnw19}.

Regarding spinors, we combine \eqref{eq:econecovdspin}, \eqref{eq:Fclp}
and the fact that $\covdg$ commutes with $F_\pm$ producing the formula
\begin{align}
  \cone{\covdg_X} (F_\pm (\cone{\Psi}))
    & = F_\pm (\covdg_X \cone{\Psi}
      \mp \tfrac{1}{2} \sqrt{\eps} \, X \clp \cone{\Psi}).
\end{align}
In other words, the pull-back of\/ $\cone{\covdg}$ along $F_\pm$ is the
Killing connection $\covda$ from \eqref{eq:covda} with Killing number
$a = \pm \tfrac{1}{2} \sqrt{\eps}$.
For a detailed proof of the following proposition, see \cite{bar96} or
\cite{sz16}.

\begin{prop} \label{coneks}
Let\/ $\Psi$ be a spinor field on $M$ and define a spinor field
$\cone{\Psi}_\pm$ on the $\eps$-metric cone $\cone{M}$ by
\begin{align} \label{eq:coneks}
  \cone{\Psi}_\pm & = \pi_2^* (F_\pm (\Psi)),
\end{align}
where we canonically identify the pull-back bundle
$\pi_2^* (\cone{\Spnr} M)$ with\/ $\Spnr \cone{M}$.
Then\/ $\cone{\Psi}_\pm$ is parallel with respect to\/ $\cone{\covdg}$
if and only if\/ $\Psi$ is Killing spinor with the Killing number
$a = \pm \tfrac{1}{2} \sqrt{\eps}$.
\end{prop}

Conversely, in order to associate a Killing spinor $\Psi$ with parallel
spinor field $\cone{\Psi}$ we need to be a bit careful and take into
account the relations \eqref{eq:spnr+1spc}, \eqref{eq:hspnr+1}.
For $n$ even, we have $\cone{\Spnr} M = \Spnr M$, hence we may choose
parallel spinor field $\cone{\Psi} \in \Sec(\Spnr\cone{M})$ arbitrarily
and it produces two Killing spinors $\Psi_\pm$, one for
each sign of the Killing number.
In the odd case we have to restrict ourselves just to half-spinor fields
$\cone{\Psi} \in \Sec(\Spnr^\pm \cone{M})$, such that each produces just
one Killing spinor $\Psi_+$ or $\Psi_-$ depending on a half-spinor
subbundle chosen.
However, this does not produce any restriction since
$\cone{\covdg}$ preserves the splitting \eqref{eq:spnr+1spc} and so we can
decompose any parallel spinor field into parallel half-spinor fields.
In any case, based on \eqref{eq:spnr+1inv} we get the formula
\begin{align} \label{eq:coneksinv}
  \Psi_\pm & = \tfrac{1}{2} \, F_\mp (\cone{\Psi} |_M).
\end{align}

Analogous equivalence for special Killing spinor-valued forms follows as
a straightforward combination of the previous two cases, see \cite{sz16}
for details.

\begin{prop} \label{conesksf}
Let\/ $\Phi$ and\/ $\Xi$ be spinor-valued differential forms on $M$ of
degrees $p$ and $p+1$ respectively, and define a spinor-valued
differential form\/ $\cone{\Theta}_\pm$ on the $\eps$-metric cone
$\cone{M}$ by
\begin{align} \label{eq:conesksf}
  \cone{\Theta}_\pm
    & = r^p \dr \wedge \pi_2^* (F_\pm (\Phi))
      + r^{p+1} \pi_2^* (F_\pm (\Xi)).
\end{align}
Then\/ $\cone{\Theta}_\pm$ is parallel with respect to\/ $\cone{\covdg}$
if and only if\/ $\Phi$ is special Killing spinor-valued form with the
corresponding $(p+1)$-form\/ $\Xi$, the Killing number
$a = \pm \tfrac{1}{2} \sqrt{\eps}$
and the Killing constant $c = \eps$.
\end{prop}

For the converse, the same considerations as in the case of spinors
apply.
Hence all special Killing spinor-valued $p$-forms with
$a = \pm \tfrac{1}{2} \sqrt{\eps}$ and $c = \eps$ are given by the
formula
\begin{align} \label{eq:conesksfinv}
  \Phi_\pm & = \tfrac{1}{2} \, F_\mp (\pi^\Norm (\cone{\Theta} |_M)),
  & \Xi_\pm & = \tfrac{1}{2} \, F_\mp (\pi^\Tan (\cone{\Theta} |_M)),
\end{align}
for all parallel spinor-valued differential forms
$\cone\Theta \in \Df^{p+1}(\cone{M}, \Spnr \cone{M})$
if $n$ is even, and
$\cone\Theta \in \Df^{p+1}(\cone{M}, \Spnr^\pm \cone{M})$
if $n$ is odd.

\section{Spaces of constant curvature}
\label{sec:ccurv}

In this section we describe the full sets of solutions of the Killing
equations of our interest on (pseudo\nh) Riemannian spaces of nonzero
constant curvature.
Without loss of generality we may assume that the sectional curvature is
equal to $\eps = \pm 1$.
We will explicitly realize the space as a quadratic hypersurface in
(pseudo\nh) Euclidean space.

Let $(n_+, n_-)$ be arbitrary signature such that $n = n_+ + n_- \ge 2$.
We take the (pseudo\nh) Euclidean space $\er^{n+1}$ with the inner
product $\cone{g}$ of signature $(\cone{n}_+, \cone{n}_-)$ given by
\eqref{eq:conesig} according to the sign of $\eps$,
\begin{align}
  \cone{g}
    & = \eps (\dif x_0)^2 + \tsum_{i=1}^{n_+} (\dif x_i)^2
        - \tsum_{j=n_+ + 1}^n (\dif x_j)^2,
\end{align}
where $x_0,\dots,x_n$ are the standard coordinates on $\er^{n+1}$.
We define manifold $\Me$ to be the connected component of the quadric
\begin{align}
  \cone{g} (x, x) = \eps ,
\end{align}
which contains the point $(1, 0, \dots, 0)$.
The manifold $\Me$ obviously inherits a (pseudo\nh) Riemannian metric
$g$ of signature $(n_+,n_-)$.

Conversely, the $\eps$-metric cone $(\cone{\Me}, \cone{g})$ over
$(M, g)$ is just a connected open submanifold of
$(\er^{n+1}, \cone{g})$, in particular, the radial coordinate is given
by
\begin{align}
  r (x) & = \sqrt{\eps \cone{g} (x, x)},
  & & \text{for all } x \in \cone{\Me}.
\end{align}
We can also identify the outer unit normal of $\Me$ with the position
vector,
\begin{align} \label{eq:ccnormal}
  \partial_r (x) & = x,
  & & \text{for all } x \in \Me.
\end{align}
The Levi-Civita connection on the cone $(\cone{M}, \cone{g})$ is simply
restriction of the ordinary partial derivative
$\cone{\covdg} = \partial$ on $\er^{n+1}$.
Using the formula \eqref{eq:conecovdg} for $\covdg$ we can verify that
$\Me$ has indeed constant sectional curvature equal to $\eps$,
\begin{align} \label{eq:cccurv}
  \curvg_{X,Y} & = \eps \rho (X \wedge Y),
  & & \text{for all } X, Y \in \Vecf(\Me),
\end{align}
where $\rho$ is the natural isomorphism from \eqref{eq:sola}.

\subsection{Killing differential forms and spinors}

The cone correspondences discussed in Propositions \ref{coneskf}, \ref{coneks},
\ref{conesksf} yield all special Killing-Yano forms,
Killing spinors and special Killing spinor-valued forms on $\Me$ with
$a = \pm \tfrac{1}{2} \sqrt{\eps}$ and $c = \eps$, by means of forms and
spinors over $\er^{n+1}$ regarded as constant sections on the cone
$\cone{\Me}$.
Substituting \eqref{eq:conefproj}, \eqref{eq:F} and \eqref{eq:ccnormal}
into \eqref{eq:coneskfinv}, \eqref{eq:coneksinv}, \eqref{eq:conesksfinv}
we get explicit formulas for the solutions,
\begin{align} \label{eq:ccskf}
\begin{split}
  \alpha (x)
    & = x \intp \cone{\theta},
\\
  \beta (x)
    & = \cone{\theta} - \dif r \wedge (x \intp \cone{\theta}),
\end{split}
  & & \text{where } \cone{\theta} \in \Wedge^{p+1} (\er^{n+1})^*,
\\ \label{eq:ccks}
  \Psi_\pm (x)
    & = \tfrac{1}{2} (1 \pm \sqrt{\eps} \, x) \clp \cone{\Psi},
  & & \text{where } \cone{\Psi} \in
    \begin{cases}
      \Spnr_{n+1}, & \text{for $n$ even,}
    \\
      \Spnr_{n+1}^\pm, & \text{for $n$ odd,}
    \end{cases}
\\ \label{eq:ccsksf}
\begin{split}
  \Phi_\pm (x)
    & = \tfrac{1}{2} (1 \pm \sqrt{\eps} \, x)
      \clp (x \intp \cone{\Theta}),
\\
  \Xi_\pm (x)
    & = \tfrac{1}{2} (1 \pm \sqrt{\eps} \, x)
      \clp (\cone{\Theta} - \dif r \wedge (x \intp \cone{\Theta})),
\\
  & \text{where } \cone{\Theta} \in
    \begin{cases}
      \Wedge^{p+1} (\er^{n+1})^* \otimes \Spnr_{n+1},
        & \text{for $n$ even,}
    \\
      \Wedge^{p+1} (\er^{n+1})^* \otimes \Spnr_{n+1}^\pm,
        & \text{for $n$ odd,}
    \end{cases}
\end{split}
  \span \span
\end{align}
for all form degrees $p = 0, \dots, n$.

As for general (not necessarily special) Killing forms, recall that
they are meaningful only in degrees $p \ge 1$.
The dimension of the solution spaces attains its universal upper bound
for Killing forms and Killing spinors by Corollaries \ref{kvfprolconn}
and \ref{ksprol}, respectively.
Hence, since $\Me$ is connected, the above formulas give all
Killing-Yano forms, Killing spinors and Killing spinor-valued forms on
$\Me$ with $a = \pm \tfrac{1}{2} \sqrt{\eps}$.
As for Killing spinors, the other Killing numbers cannot occur
by~\eqref{eq:ksscal}.
Hence it remains to discuss Killing spinor-valued forms with
$a \ne \pm \tfrac{1}{2} \sqrt{\eps}$.

\subsection{Other Killing numbers}

To determine possible Killing numbers $a$ admitting nontrivial Killing
spinor-valued forms is more involved.
Relying on the first integrability condition \eqref{eq:kvfint0},
we employ separately the curvature $\curvaff = \curvg$ acting
on the form part and $\curvv = \curva$ acting on the value (spinor)
part.
By \eqref{eq:cccurv} and \eqref{eq:solaform} we have
\begin{align}
  \curvaff_{X, Y} \Phi
    & = \eps (X^\flat \wedge (Y \intp \Phi)
            - Y^\flat \wedge (X \intp \Phi)),
\end{align}
and using \eqref{eq:skewcurv1}, \eqref{eq:skewcurv2} and
\eqref{eq:curvdecomp3} we compute
\begin{align}
\begin{split}
  \curvaff_X \wedge \Phi
    & = \eps \tsum_{j=1}^n e^j
        \wedge (X^\flat \wedge (e_j \intp \Phi)
                - (e_j)^\flat \wedge (X \intp \Phi))
\\
    & = \eps \Big(
        -X^\flat \wedge \tsum_{j=1}^n e^j \wedge (e_j \intp \Phi)
        - \tsum_{j,k=1}^n g_{jk} \, e^j \wedge e^k \wedge (X \intp \Phi)
      \Big)
\\
    & = - \eps p \, X^\flat \wedge \Phi,
\end{split}
\\
  \curvaff \wedge \Phi
    & = - \eps p \tsum_{i=1}^n e^i \wedge (e_i)^\flat \wedge \Phi
    = - \eps p \tsum_{i,j=1}^n g_{ij} \, e^i \wedge e^j \wedge \Phi
    = 0,
\\ \label{eq:ccp2curvaff}
\begin{split}
  (\curvaff \Phi)^{(p,2)}_{X,Y}
    & = \eps (
        X^\flat \wedge (Y \intp \Phi)
        - Y^\flat \wedge (X \intp \Phi)
        + {}
\\
    & \qquad
        + Y \intp (X^\flat \wedge \Phi)
        - X \intp (Y^\flat \wedge \Phi)
      )
    = 0.
\end{split}
\end{align}
For the value part we have by \eqref{eq:curva}, \eqref{eq:cccurv} and
\eqref{eq:spinla},
\begin{align} \label{eq:cccurvv}
  \curvv_{X,Y} \Phi
    & = -\tfrac{1}{2} (\eps - 4 a^2) (X \clp {Y \clp} + g (X, Y)) \, \Phi,
\end{align}
and using again \eqref{eq:skewcurv1}, \eqref{eq:skewcurv2},
\eqref{eq:curvdecomp3} and also \eqref{eq:clf} we compute
\begin{align}
\begin{split}
  \curvv_X \wedge \Phi
    & = -\tfrac{1}{2} (\eps - 4 a^2)
        \tsum_{j=1}^n e^j \wedge (
            (X \clp {e_j \clp} + g (X, e_j)) \, \Phi
          )
\\
    & = -\tfrac{1}{2} (\eps - 4 a^2)
        (X \clp (\clf \wedge \Phi) + X^\flat \wedge \Phi),
\end{split}
\\
\begin{split}
  \curvv \wedge \Phi
    & = -\tfrac{1}{4} (\eps - 4 a^2)
        \tsum_{i=1}^n e^i \wedge (e_i \clp (\clf \wedge \Phi)
                                  + (e_i)^\flat \wedge \Phi)
\\
    & = -\tfrac{1}{4} (\eps - 4 a^2) \Big(
        \clf \wedge \clf \wedge \Phi
        + \tsum_{i,j=1}^n g_{ij} \, e^i \wedge e^j \wedge \Phi
      \Big)
\\
    & = -\tfrac{1}{4} (\eps - 4 a^2) \, \clf \wedge \clf \wedge \Phi,
\end{split}
\\ \label{eq:ccp2curvv}
\begin{split}
  (\curvv \Phi)^{(p,2)}_{X,Y}
    & = -\tfrac{1}{2} (\eps - 4 a^2) \big(
        (X \clp {Y \clp} + g (X, Y)) \, \Phi
        - {}
\\
    & \quad
        - \tfrac{1}{p} \big(
          Y \intp (X \clp (\clf \wedge \Phi) + X^\flat \wedge \Phi)
          - {}
\\
    & \quad \qquad
          - X \intp (Y \clp (\clf \wedge \Phi) + Y^\flat \wedge \Phi)
          - {}
\\
    & \quad \qquad
          - \tfrac{1}{p+1} \, Y \intp (
              X \intp (\clf \wedge \clf \wedge \Phi)
            )
        \big)
      \big).
\end{split}
\end{align}
Altogether, on $\Me$ we have
$(\curv\Phi)^{(p,2)} = (\curvv\Phi)^{(p,2)}$ by the previous formulas.

Now we shall prove that there are no nontrivial solutions with Killing
number $a \ne \pm \tfrac{1}{2} \sqrt{\eps}$ for $p \ge 2$.
Recall that for $p = 1$ the component $ (\curv\Phi)^{(p,2)}$ of the
curvature action vanishes automatically and as we shall observe
later on,
there exist additional Killing spinor-valued 1-forms on $\Me$.

\begin{lem} \label{ccksfint01}
Let\/ $\Phi$ be a Killing spinor-valued $p$-form on $\Me$ with
$a \ne \pm \tfrac{1}{2} \sqrt{\eps}$ and $p \ge 2$.
Then it holds
\begin{align}
  \clv \intp (\clv \intp \Phi) & = 0.
\end{align}
\end{lem}

\begin{prf}
We compute the following curvature operator built from
$(\curv\Phi)^{(p,2)}$,
\begin{align*}
  & r_1 ((\curv\Phi)^{(p,2)})
    = \tsum_{i,j,k,l=1}^n g^{ik} g^{jl} \,
          e_i \intp \big( e_j \intp (\curv\Phi)^{(p,2)}_{e_k,e_l} \big)
\\
  & \quad
    = -\tfrac{1}{2} (\eps - 4 a^2) \Big(
        \tsum_{i,j,k,l=1}^n g^{ik} g^{jl} \, e_i \intp (e_k \clp (
          e_j \intp ( e_l \clp \Phi)
        ))
        + {}
\\
  & \qquad \qquad \qquad \qquad
        + \tsum_{i,j=1}^n g^{ij} \, e_i \intp (e_j \intp \Phi)
      \Big)
\\
  & \quad
    =  -\tfrac{1}{2} (\eps - 4 a^2) \, \clv \intp (\clv \intp \Phi).
\end{align*}
Now the claim follows from Proposition \ref{kvfprol}.
\end{prf}

\begin{lem} \label{ccksfint02}
Let\/ $\Phi$ be a Killing spinor-valued $p$-form on $\Me$ with
$a \ne \pm \tfrac{1}{2} \sqrt{\eps}$ and $p \ge 2$.
Then it holds
\begin{align}
  \clv \intp \Phi & = 0.
\end{align}
\end{lem}

\begin{prf}
Again we compute a curvature operator built from $(\curv\Phi)^{(p,2)}$,
\begin{align*}
  & r_2 ((\curv\Phi)^{(p,2)})
    = \tsum_{i,j,k,l=1}^n g^{ik} g^{jl} \,
          e_i \clp \big( e_j \intp (\curv\Phi)^{(p,2)}_{e_k,e_l} \big)
\\
  & \quad
    = -\tfrac{1}{2} (\eps - 4 a^2) \Big(
        \tsum_{i,j,k,l=1}^n g^{ik} g^{jl} \, e_i \clp e_k \clp (
          e_j \intp (e_l \clp \Phi)
        )
        + {}
\\
  & \qquad \qquad \qquad \qquad
        + \tsum_{i,j=1}^n g^{ij} \, e_i \clp (e_j \intp \Phi)
        - {}
\\
  & \qquad \qquad \qquad \qquad
        - \tfrac{1}{p} \Big(
          \tsum_{i,j,k,l=1}^n g^{ik} g^{jl} \, e_k \intp (e_i \clp (
            e_j \intp (e_l \clp (\clf \wedge \Phi))
          ))
          + {}
\\
  & \qquad \qquad \qquad \qquad \qquad \;
          + \tsum_{i,j,k,l=1}^n g^{ik} g^{jl} \, e_k \intp (e_i \clp (
            e_j \intp ((e_l)^\flat \wedge \Phi)
          ))
        \Big)
      \Big)
\\
  & \quad
    = \tfrac{1}{2 p} (\eps - 4 a^2) (
        (n p + n - 2 p) \, \clv \intp \Phi
        + \clv \intp (\clv \intp (\clf \wedge \Phi))
      ),
\\
\intertext{and rearrange the second term using \eqref{eq:clfclv},}
  & \quad
    = \tfrac{1}{2 p} (\eps - 4 a^2) (
        (n + 2) (p - 1)\, \clv \intp \Phi
        + \clf \wedge (\clv \intp (\clv \intp \Phi))
      ).
\end{align*}
Now the claim follows from Proposition \ref{kvfprol} and Lemma
\ref{ccksfint01}.
\end{prf}

\begin{prop} \label{ccksfint0}
There are no nontrivial Killing spinor-valued $p$-forms on $\Me$ with
$a \ne \pm \tfrac{1}{2} \sqrt{\eps}$ and $p \ge 2$.
\end{prop}

\begin{prf}
Again we compute a curvature operator built from $(\curv\Phi)^{(p,2)}$,
\begin{align*}
  & r_3 ((\curv\Phi)^{(p,2)})
    = \tsum_{i,j,k,l=1}^n g^{ik} g^{jl} \,
          e_i \clp e_j \clp (\curv\Phi)^{(p,2)}_{e_k,e_l}
\\
  & \quad
    = -\tfrac{1}{2} (\eps - 4 a^2) \Big(
        - \tsum_{i,j,k,l=1}^n g^{ik} g^{jl} \,
            e_i \clp e_k \clp e_j \clp e_l \clp \Phi
        - {}
\\
  & \qquad \qquad \qquad
        - 2 \tsum_{i,l=1}^n g^{il} \, e_i \clp e_l \clp \Phi
        + {}
\\
  & \qquad \qquad
        + \tsum_{i,j=1}^n g^{ij} \, e_i \clp e_j \clp \Phi
        - {}
\\
  & \qquad
        - \tfrac{1}{p} \Big(
          - \tsum_{i,j,k,l=1}^n g^{ik} g^{jl} \,
            e_i \clp e_k \clp (
              e_l \intp (e_j \clp (\clf \wedge \Phi))
            )
          - {}
\\
  & \qquad \qquad \; \,
          - 2 \tsum_{i,l=1}^n g^{il} \,
            e_l \intp (e_i \clp (\clf \wedge \Phi))
          - {}
\\
  & \qquad \qquad \; \,
          - \tsum_{i,j,k,l=1}^n g^{ik} g^{jl} \,
            (e_k)^\flat \wedge (e_i \clp (
              e_l \intp (e_j \clp \Phi)
            ))
          + {}
\\
  & \qquad \qquad \; \,
          + \tsum_{i,j=1}^n g^{ij} \, e_i \clp e_j \clp \Phi
          - {}
\\
  & \qquad \qquad \; \,
          - \tsum_{i,j,k,l=1}^n g^{ik} g^{jl} \,
            e_k \intp (e_i \clp
              e_j \clp e_l \clp (\clf \wedge \Phi)
            )
          - {}
\\
  & \qquad \qquad \; \,
          - \tsum_{i,j,k,l=1}^n g^{ik} g^{jl} \,
            e_k \intp (e_i \clp (
              (e_l)^\flat \wedge (e_j \clp \Phi)
            ))
          + {}
\\
  & \qquad \qquad \; \,
          + \tfrac{1}{p+1} \tsum_{i,j,k,l=1}^n g^{ik} g^{jl} \,
            e_k \intp (e_i \clp (
              e_l \intp (e_j \clp (\clf \wedge \clf \wedge \Phi))
            ))
        \Big)
      \Big)
\\
  & \quad
    = \tfrac{1}{2 p} (\eps - 4 a^2) \big(
        n (n p - p - 1) \, \Phi
        - \clf \wedge (\clv \intp \Phi)
        + {}
\\
  & \qquad \qquad \qquad \quad \; \,
        + (2 n - 3) \, \clv \intp (\clf \wedge \Phi)
        + {}
\\
  & \qquad \qquad \qquad \quad \; \,
        + \tfrac{1}{p+1} \, \clv \intp (\clv \intp (
            \clf \wedge \clf \wedge \Phi
          ))
      \big),
\\
\intertext{and rearrange the last two terms using \eqref{eq:clfclv},}
  & \quad
    = \tfrac{1}{2 (p+1)} (\eps - 4 a^2) \big(
        (n + 1) (n + 2) (p - 1) \, \Phi
        + {}
\\
  & \qquad \qquad \qquad \qquad \quad
        + \tfrac{1}{p} (
          2 (n + 2) (p - 1) \, \clf \wedge (\clv \intp \Phi)
          + {}
\\
  & \qquad \qquad \qquad \qquad \qquad \quad
          + \clf \wedge \clf \wedge (\clv \intp (\clv \intp \Phi))
        )
      \big).
\end{align*}
Now the claim follows from Proposition \ref{kvfprol} and Lemma
\ref{ccksfint02}.
\end{prf}

Before we proceed to the degree 1 case we conclude this part discussing
all the possible special Killing spinor-valued forms in degree 0 on
$\Me$.

\begin{prop} \label{ccsksf0}
Let\/ $\Phi$ be a nontrivial special Killing spinor-valued 0-form on
$\Me$ and\/ $\Xi$ the corresponding spinor-valued 1-form.
Then the Killing number is necessarily
$a = \pm \tfrac{1}{2} \sqrt{\eps}$, and the constant $c$ is either
\begin{enumerate}
\item[a)]
  $c = \eps$, in which case\/ $\Phi$ and\/ $\Xi$ are given by the
  formula \eqref{eq:ccsksf}, or
\item[b)]
  $c = 0$, in which case\/ $\Phi$ is a Killing spinor and\/ $\Xi = 0$.
\end{enumerate}
\end{prop}

\begin{prf}
First we employ the first integrability condition \eqref{eq:skvfint0}
in Proposition \ref{skvfint0}.
Since $\Phi$ is of degree 0, the right hand side vanishes and by
\eqref{eq:cccurvv} the equation \eqref{eq:skvfint0} reads
\begin{align*}
  (\eps - 4 a^2) \, \rho (X \wedge Y) \, \Phi
    = - \tfrac{1}{2} (\eps - 4 a^2) \,
      (X \clp {Y \clp} + g (X, Y)) \, \Phi & = 0.
\end{align*}
Note that this is the same condition as \eqref{eq:ksint0} for
Killing spinors and we must again have
$a = \pm \tfrac{1}{2} \sqrt{\eps}$.
In more detail, we can argue that the spin representation of the spin
Lie algebra contains no trivial summands.
Alternatively we can just compute the operator $r_3$ as in Proposition
\ref{ccksfint0},
\begin{align*}
  0 = r_3(\curv \Phi)
    & = \tfrac{1}{2} (\eps - 4 a^2) n (n - 1) \, \Phi.
\end{align*}

As for the second part we employ the second integrability condition
\eqref{eq:skvfint1} in Proposition \ref{skvfint1}.
Because $a = \pm \tfrac{1}{2} \sqrt{\eps}$
we have also $\curvv = \curva = 0$ and by \eqref{eq:cccurv} the equation
\eqref{eq:skvfint1} reads
\begin{align*}
  (\eps - c) \, \rho(X \wedge Y)^\aff \, \Xi
    = (\eps - c) (X \intp (Y^\flat \wedge \Xi)
                  - Y \intp (X^\flat \wedge \Xi))
    & = 0.
\end{align*}
Hence we must have either $c = \eps$ or $\Xi = 0$.
Again we can argue that the representation of the spin Lie algebra on
spinor-valued 1-forms contains no trivial summands, or to compute the
operator $q$ (sometimes called the \emph{curved Casimir~operator}),
\begin{align*}
\begin{split}
  0 & = q((\curv - c \rho^\aff) \, \Xi)
    = (\eps - c) \, q (\rho^\aff \Xi)
\\
    & = (\eps - c) \tsum_{i,j,k,l=1}^n g^{ik} g^{il} \,
      (e_i)^\flat \wedge (e_j \intp (
        \rho^\aff (e_k \wedge e_l) \, \Xi
      ))
\\
    & = - (\eps - c) (n - 1) \, \Xi.
\end{split}
\end{align*}
Finally, the case $\Xi = 0$ implies $c = 0$ by the second defining
equation \eqref{eq:skvf}.
\end{prf}

\subsection{Additional solutions in degree 1}

To resolve the case of Killing spinor-valued 1-forms on the space $\Me$
of constant curvature we need to employ also the second integrability
condition \eqref{eq:kvfint1}.
So let $\Phi'$ be a Killing spinor-valued 1-form with Killing number
$a' \ne \pm \tfrac{1}{2} \sqrt{\eps}$ and the corresponding spinor-valued
2-form $\Xi'$.
The left hand side of \eqref{eq:kvfint1} is just a multiple of the
$(2,2)$-symmetry type component, so by \eqref{eq:vfcurvmod},
\eqref{eq:ccp2curvaff} and \eqref{eq:ccp2curvv} we~have
\begin{align} \label{eq:cc22mcurv}
\begin{split}
  (\mcurv \Xi')^{(2,2)}_{X,Y}
    = -\tfrac{3}{2} (\eps - 4 (a')^2) \big( &
        (X \clp {Y \clp} + g (X, Y)) \, \Phi'
        - {}
\\
        & - \tfrac{1}{2} \big(
          Y \intp (X \clp (\clf \wedge \Phi') + X^\flat \wedge \Phi')
          - {}
\\
        & \qquad
          - X \intp (Y \clp (\clf \wedge \Phi') + Y^\flat \wedge \Phi')
          - {}
\\
        & \qquad
          - \tfrac{1}{3} \, Y \intp (
              X \intp (\clf \wedge \clf \wedge \Phi')
            )
        \big)
      \big),
\end{split}
\end{align}
noting that the degree of $\Xi'$ is $p + 1 = 2$.
In order to compute the right hand side we need the covariant derivative
of the curvature.
Note that the isomorphism $\rho$ from \eqref{eq:sola} and
\eqref{eq:spinla} is invariant with respect to the spin group and hence
we have $\covdg \rho = 0$.
By \eqref{eq:cccurv} and \eqref{eq:cccurvv} we have
\begin{align}
  (\covd_X \curvaff)_{Y,Z} & = \eps \, (\covdg_X \rho) (Y \wedge Z) = 0
\\
\begin{split}
  (\covd_X \curvv)_{Y,Z}
    & = (\eps - 4 (a')^2) (\covdap_X \rho) (Y \wedge Z)
\\
    & = (\eps - 4 (a')^2) \big(
        (\covdg_X \rho) (Y \wedge Z)
        - a'[{X \clp}, \rho (Y \wedge Z)]
      \big)
\\
    & = a' (\eps - 4 (a')^2)
        {(g (X, Z) \, Y - g (X, Y) \, Z) \clp},
\end{split}
\end{align}
for all $X, Y, Z \in \Vecf(\Me)$.
Next we compute the particular action of $\covd \curv$ on $\Phi'$ defined
by the equation \eqref{eq:barwedge},
\begin{align}
\begin{split}
  (\covd_X \curv)_Y \wedge \Phi'
    & = a' (\eps - 4 (a')^2) \tsum_{j=1}^n
      e^j \wedge (g (X, e_j) \, Y - g (X, Y) \, e_j) \clp \Phi'
\\
    & = a' (\eps - 4 (a')^2) (X^\flat \wedge (Y \clp \Phi')
                        - g (X, Y) \, \clf \wedge \Phi'),
\end{split}
\\
\begin{split}
  (\covd_X \curv) \wedge \Phi'
    & = \tfrac{1}{2} \, a' (\eps - 4 (a')^2) \tsum_{i=1}^n
      e^i \wedge (X^\flat \wedge (e_i \clp \Phi')
                  - {}
\\
  & \qquad \qquad \qquad \qquad \qquad \qquad
                  - g (X, e_i) \, \clf \wedge \Phi'),
\\
    & = - a' (\eps - 4 (a')^2) \, X^\flat \wedge (\clf \wedge \Phi'),
\end{split}
\\ \label{eq:ccbarwedge}
\begin{split}
  ((\covd \curv) \barwedge \Phi')_{X, Y}
    & = a' (\eps - 4 (a')^2) \big(
        X^\flat \wedge (Y \clp \Phi') - Y^\flat \wedge (X \clp \Phi')
        + {}
\\
    & \qquad \qquad \qquad \quad \;
        + \tfrac{1}{2} (
          Y \intp (X^\flat \wedge \clf \wedge \Phi')
          - {}
\\
    & \qquad \qquad \qquad \qquad \quad \;
          - X \intp (Y^\flat \wedge \clf \wedge \Phi')
        )
      \big).
\end{split}
\end{align}
Now we proceed similarly to Lemmas \ref{ccksfint01}, \ref{ccksfint02}
and Proposition \ref{ccksfint0} and compare the operators $r1$, $r2$ and
$r3$ applied to both sides of \eqref{eq:kvfint1}.

\begin{lem} \label{ccksfint11}
Let\/ $\Phi'$ be a Killing spinor-valued 1-form on $\Me$ with
$a' \ne \pm \tfrac{1}{2} \sqrt{\eps}$ and\/ $\Xi'$ the corresponding
spinor-valued 2-form.
Then it holds
\begin{align}
  \clv \intp (\clv \intp \Xi')
    & = \tfrac{4}{3} \, a' (n-1) \, \clv \intp \Phi'.
\end{align}
\end{lem}

\begin{prf}
From \eqref{eq:cc22mcurv} and the computation in Lemma \ref{ccksfint01}
we have
\begin{align*}
  r_1 ((\mcurv \Xi')^{(2,2)})
    & = - \tfrac{3}{2} (\eps - 4 (a')^2) \, \clv \intp (\clv \intp \Xi').
\end{align*}
For the right hand side we compute using \eqref{eq:ccbarwedge},
\begin{align*}
  & r_1 ((\covd \curv) \barwedge \Phi') = {}
\\
  & \quad
    = a' (\eps - 4 (a')^2) \Big(
        - \tsum_{i,j,k,l=1}^n g^{ik} g^{jl} \,
          e_i \intp ((e_k)^\flat \wedge (
            e_j \intp (e_l \clp \Phi')
          ))
        + {}
\\
  & \qquad \qquad \qquad \qquad \; \; \,
        + \tsum_{i,l=1}^n g^{il} \,
          e_i \intp (e_l \clp \Phi')
        - {}
\\
  & \qquad \qquad \qquad \qquad \; \; \,
        - \tsum_{i,j,k,l=1}^n g^{ik} g^{jl} \,
          e_i \intp (e_k \clp (
            e_j \intp ((e_l)^\flat \wedge \Phi')
          ))
      \Big)
\\
  & \quad
    = - 2 a' (\eps - 4 (a')^2) (n - 1) \, \clv \intp \Phi'.
\end{align*}
Now the claim follows from Proposition \ref{kvfint1}.
\end{prf}

\begin{lem} \label{ccksfint12}
Let\/ $\Phi'$ be a Killing spinor-valued 1-form on $\Me$ with
$a' \ne \pm \tfrac{1}{2} \sqrt{\eps}$
and\/ $\Xi'$ the corresponding spinor-valued 2-form.
Then it holds
\begin{align}
  \clv \intp \Xi'
    & = \tfrac{2}{3} \, a' ((n-2) \, \Phi'
        - \clf \wedge (\clv \intp \Phi')).
\end{align}
\end{lem}

\begin{prf}
From \eqref{eq:cc22mcurv} and the computation in Lemma \ref{ccksfint02}
we have
\begin{align*}
  r_2 ((\mcurv \Xi')^{(2,2)})
    & = \tfrac{3}{4} (\eps - 4 (a')^2) (
        (n + 2) \, \clv \intp \Xi'
        + \clf \wedge (\clv \intp (\clv \intp \Xi'))
      ).
\end{align*}
For the right hand side we compute using \eqref{eq:ccbarwedge},
\begin{align*}
  & r_2 ((\covd \curv) \barwedge \Phi') = {}
\\
  & \quad
    = a' (\eps - 4 (a')^2) \Big(
        - \tsum_{i,j,k,l=1}^n g^{ik} g^{jl} \,
          (e_k)^\flat \wedge (e_i \clp (
            e_j \intp (e_l \clp \Phi')
          ))
        + {}
\\
  & \qquad \qquad \qquad \qquad \; \; \,
        + \tsum_{i,l=1}^n g^{il} \,
          e_i \clp e_l \clp \Phi'
        - {}
\\
  & \qquad \qquad \qquad \qquad \; \; \,
        - \tsum_{i,j,k,l=1}^n g^{ik} g^{jl} \,
          e_i \clp e_k \clp (
            e_j \intp ((e_l)^\flat \wedge \Phi')
          )
        + {}
\\
  & \qquad \qquad \qquad \qquad \; \; \,
        + \tfrac{1}{2} \tsum_{i,j,k,l=1}^n g^{ik} g^{jl} \,
          e_k \intp (e_i \clp (
            e_j \intp ((e_l)^\flat \wedge \clf \wedge \Phi')
          ))
      \Big)
\\
  & \quad
    = a' (\eps - 4 (a')^2) \big(
        n (n - 2) \, \Phi'
        - \clf \wedge (\clv \intp \Phi')
        + \tfrac{1}{2} (n - 2) \, \clv \intp (\clf \wedge \Phi')
      \big),
\\
\intertext{and rearrange the last term using \eqref{eq:clfclv},}
  & \quad
    = \tfrac{1}{2} \, a' (\eps - 4 (a')^2) (
        (n - 2) (n + 2) \, \Phi'
        + (n - 4) \, \clf \wedge (\clv \intp \Phi')
      ).
\end{align*}
Now the claim follows from Proposition \ref{kvfint1} and Lemma
\ref{ccksfint11}.
\end{prf}

\begin{prop} \label{ccksfint1}
Let\/ $\Phi'$ be a Killing spinor-valued 1-form on $\Me$ with
$a' \ne \pm \tfrac{1}{2} \sqrt{\eps}$
and\/ $\Xi'$ the corresponding spinor-valued 2-form.
Then it holds
\begin{align}
  \Xi' & = - \tfrac{2}{3} \, a' \clf \wedge \Phi'.
\end{align}
In other words, $\Phi'$ satisfies differential equation
\begin{align} \label{eq:32ksf}
  \covdg_X \Phi' & = a' \big(
        X \clp \Phi' - \tfrac{2}{3} \, X \intp (\clf \wedge \Phi')
      \big),
  & & \text{for all } X \in \Vecf(\Me).
\end{align}
\end{prop}

\begin{prf}
From \eqref{eq:cc22mcurv} and the computation in Proposition
\ref{ccksfint0} we have
\begin{align*}
\begin{split}
  r_3 ((\mcurv \Xi')^{(2,2)})
    & = \tfrac{1}{2} (\eps - 4 (a')^2) \big(
        (n + 1) (n + 2) \, \Xi'
        + (n + 2) \, \clf \wedge (\clv \intp \Xi')
        + {}
\\ & \qquad \qquad \qquad \quad \,
        + \tfrac{1}{2} \, \clf \wedge \clf
            \wedge (\clv \intp (\clv \intp \Xi'))
      \big).
\end{split}
\end{align*}
For the right hand side we compute using \eqref{eq:ccbarwedge},
\begin{align*}
  & r_3 ((\covd \curv) \barwedge \Phi') = {}
\\
  & \quad
    = a' (\eps - 4 (a')^2) \Big(
        \tsum_{i,j,k,l=1}^n g^{ik} g^{jl} \,
          (e_k)^\flat \wedge (e_i \clp
            e_j \clp e_l \clp \Phi'
          )
        + {}
\\
  & \qquad \qquad \qquad \qquad \; \; \,
        + \tsum_{i,j,k,l=1}^n g^{ik} g^{jl} \,
          e_i \clp e_k \clp (
            (e_l)^\flat \wedge (e_j \clp \Phi')
          )
        + {}
\\
  & \qquad \qquad \qquad \qquad \; \; \,
        + 2 \tsum_{i,l=1}^n g^{il} \,
          (e_l)^\flat \wedge (e_i \clp \Phi')
        + {}
\\
  & \qquad \qquad \qquad \qquad \; \; \,
        + \tfrac{1}{2} \Big(
          - \tsum_{i,j,k,l=1}^n g^{ik} g^{jl} \,
            (e_k)^\flat \wedge (e_i \clp (
              e_l \intp (e_j \clp (\clf \wedge \Phi'))
            ))
          + {}
\\
  & \qquad \qquad \qquad \qquad \qquad \quad \,
          + \tsum_{i,j=1}^n g^{ij} \,
              e_i \clp e_j \clp (\clf \wedge \Phi')
          - {}
\\
  & \qquad \qquad \qquad \qquad \qquad \quad \,
          - \tsum_{i,j,k,l=1}^n g^{ik} g^{jl} \,
            e_k \intp (e_i \clp (
              (e_l)^\flat \wedge (e_j \clp (\clf \wedge \Phi'))
            ))
        \Big)
      \Big)
\\ & \quad
    = - \tfrac{1}{2} \, a' (\eps - 4 (a')^2) (
        (5 n - 4) \, \clf \wedge \Phi'
        + {}
\\ & \qquad \qquad \qquad \qquad \qquad
        + \clf \wedge (\clv \intp (\clf \wedge \Phi'))
        + \clv \intp (\clf \wedge \clf \wedge \Phi')
      ),
\\
\intertext{and rearrange the last two terms using \eqref{eq:clfclv},}
  & \quad
    = - a' (\eps - 4 (a')^2) (
      (n + 2) \, \clf \wedge \Phi'
      + \clf \wedge \clf \wedge (\clv \intp \Phi')
      ).
\end{align*}
Now the claim follows from Proposition \ref{kvfint1} and Lemma
\ref{ccksfint12}.
\end{prf}

Consequently, it remains to describe solutions of the stronger equation
\eqref{eq:32ksf}.
It turns out that for $a' \ne 0$ the solutions are just algebraic
transformation of suitable special Killing spinor-valued forms in degree
0.
The transformation works in general, so we point out that the following
two propositions apply to any spin (pseudo\nh) Riemannian manifold and
not just $M_\eps$.

\begin{prop} \label{32ksfcor}
Let $M$ be an arbitrary spin (pseudo\nh) Riemannian manifold.
Then spinor-valued differential 1-forms\/ $\Phi'$ on $M$ solving
\eqref{eq:32ksf} with Killing number $a' \ne 0$ bijectively correspond
to special Killing spinor-valued 0-forms\/ $\Phi$ on $M$ with Killing
number $a = \tfrac{1}{3} \, a'$, constant $c = 4a^2$ and the
corresponding spinor-valued 1-form\/ $\Xi$ such that
\begin{align} \label{eq:32ksfprimit}
  \clv \intp \Xi & = - 2 a \Phi.
\end{align}
The correspondence is given by formulas
\begin{align} \label{eq:32ksfcor1}
  \Phi & = - \tfrac{1}{2 a (n + 1)} \, \clv \intp \Phi',
  & \Xi & = \Phi' + \tfrac{1}{n+1} \, \clf \wedge (\clv \intp \Phi'),
\\ \label{eq:32ksfcor2}
  \Phi' & = \Xi + 2 a \clf \wedge \Phi
    = \Xi - \clf \wedge (\clv \intp \Xi).
  \span \span
\end{align}
\end{prop}

\begin{prf}
Let us first note that since $\Phi'$ has degree 1 we have
$\clv \intp (\clv \intp \Phi') = 0$ and
$X \intp (\clv \intp \Phi') = 0$,
which we shall use repeatedly in the proof.
This gives an additional insight why there is no analogous construction
in higher degrees.

Now let $\Phi'$ be a solution of \eqref{eq:32ksf} with $a' \ne 0$ and
define $\Phi$ and $\Xi$ by the formulas \eqref{eq:32ksfcor1}.
Using \eqref{eq:clfclv} we immediately get the equation
\eqref{eq:32ksfprimit},
\begin{align*}
  \clv \intp \Xi
    & = \clv \intp \Phi' - \tfrac{n}{n+1} \, \clv \intp \Phi'
    = - 2 a \Phi.
\end{align*}
Then we compute covariant derivatives using the assumption
\eqref{eq:32ksf}, formulas for computing with spinor-valued forms
\eqref{eq:clf1}--\eqref{eq:clf4}, \eqref{eq:clfclv}, and also
$a = \tfrac{1}{3} \, a'$,
\begin{align*}
\begin{split}
  \covdg_X (\clv \intp \Phi')
    & = \clv \intp (\covdg_X \Phi')
\\
    & = a' \big(
        - 2 X \intp \Phi' - X \clp (\clv \intp \Phi')
        + \tfrac{2}{3} \, X \intp (\clv \intp (\clf \wedge \Phi'))
      \big)
\\
    & = - a (
        X \clp (\clv \intp \Phi')
        + 2 (n + 1) \, X \intp \Phi'
      ),
\end{split}
\\
\begin{split}
  \covdg_X (\clf \wedge (\clv \intp \Phi'))
    & = \clf \wedge (\covdg_X (\clv \intp \Phi'))
\\
    & = - a (
        -2 X^\flat \wedge (\clv \intp \Phi')
        - X \clp (\clf \wedge (\clv \intp \Phi'))
        + {}
\\
  & \qquad \; \; \;
        + 2 (n + 1) (
          X \clp \Phi'
          - X \intp (\clf \wedge \Phi')
        )
      ).
\end{split}
\end{align*}
The defining equations \eqref{eq:ksf}, \eqref{eq:sksf} of special
Killing spinor-valued forms then follow using also \eqref{eq:32ksfcor1}
and $c = 4 a^2$,
\begin{align*}
\begin{split}
  \covdg_X \Phi
    & = - a X \clp \Phi + X \intp \Xi
      - \tfrac{1}{n+1} \, X \intp (\clf \wedge (\clv \intp \Phi'))
\\
    & = a X \clp \Phi + X \intp \Xi,
\end{split}
\\
\begin{split}
  \covdg_X \Xi
    & = a \big(
        3 X \clp \Phi' - 2 X \intp (\clf \wedge \Phi')
        + {}
\\
  & \qquad
        + \tfrac{2}{n+1} \, X^\flat \wedge (\clv \intp \Phi')
        + \tfrac{1}{n+1} \, X \clp (\clf \wedge (\clv \intp \Phi'))
        - {}
\\
  & \qquad
        - 2 X \clp \Phi'
        + 2 X \intp (\clf \wedge \Phi')
      \big)
\\
    & = a X \clp \Xi - c X^\flat \wedge \Phi.
\end{split}
\end{align*}

Conversely, suppose that $\Phi$ and $\Xi$ satisfy \eqref{eq:ksf},
\eqref{eq:sksf} with $a \ne 0$ and $c = 4 a^2$.
We again compute covariant derivative using the formulas
\eqref{eq:clf1}--\eqref{eq:clf4} and \eqref{eq:clfclv},
\begin{align*}
\begin{split}
  \covdg_X (\clf \wedge \Phi)
    & = \clf \wedge (\covdg_X \Phi)
\\
    & = a (- 2 X^\flat \wedge \Phi - X \clp (\clf \wedge \Phi))
      + X \clp \Xi  - X \intp (\clf \wedge \Xi).
\end{split}
\end{align*}
Now we define $\Phi'$ by the formula \eqref{eq:32ksfcor2} and further
compute using also $a' = 3 a$,
\begin{align*}
\begin{split}
  \covdg_X \Phi'
    & = a X \clp \Xi - 4 a^2 \, X^\flat \wedge \Phi
        - {}
\\
  & \quad
        - 2 a^2 (
          2 X^\flat \wedge \Phi + X \clp (\clf \wedge \Phi)
        )
        + 2 a (
          X \clp \Xi - X \intp (\clf \wedge \Xi)
        )
\\
    & = 3 a X \clp \Xi - 2 a X \intp (\clf \wedge \Xi)
        + {}
\\
  & \quad
        + 6 a^2 \, X \clp (\clf \wedge \Phi)
        - 4 a^2 \, X \intp (\clf \wedge \clf \wedge \Phi)
\\
    & = a' (X \clp \Phi'
        - \tfrac{2}{3} \, X \intp (\clf \wedge \Phi')),
\end{split}
\end{align*}
proving that $\Phi'$ solves \eqref{eq:32ksf}.
Finally, a straightforward computation using \eqref{eq:32ksfprimit} and
\eqref{eq:clfclv} verifies that the formulas \eqref{eq:32ksfcor1} and
\eqref{eq:32ksfcor2} are inverse to each other.
\end{prf}

The last proposition allows to translate the cone correspondence from
Proposition \ref{conesksf} to Killing spinor-valued 1-forms satisfying
\eqref{eq:32ksf}.
As it turns out, the condition \eqref{eq:32ksfprimit} has a nice
representation theoretic formulation in terms of the corresponding
parallel spinor-valued 1-form $\cone{\Theta}$ on the cone $\cone{M}$.
In order to see it, we recall the well-known invariant decomposition of
spinor-valued 1-forms corresponding on the level of vector spaces to
\begin{align}
  (\er^n)^* \otimes \Spnr_n
    & \simeq \Spnr_n \oplus \Spnr^\Tw_n.
\end{align}
The respective projections are given by
\begin{align}
\begin{split}
  \pi^\Spnr (\Theta)
    & = \clv \intp \Theta,
\\
  \pi^\Tw (\Theta)
    & = \Theta + \tfrac{1}{n} \, \clf \wedge (\clv \intp \Theta),
\end{split}
  & & \text{for all } \Theta \in (\er^n)^* \otimes \Spnr_n,
\end{align}
and the space $\Spnr^\Tw$ of \emph{primitive} spinor-valued 1-forms is
simply the kernel of $\pi^\Spnr$.
In even dimensions we define also the spaces $\Spnr^{\Tw\pm}$ of
primitive half-spinor-valued 1-forms as the kernel of $\pi^\Spnr$
restricted to half-spinor-valued forms.

\begin{prop} \label{cone32ksf}
Let $M$ be a spin (pseudo\nh) Riemannian manifold and\/
$\Phi'$ be a spinor-valued differential 1-form on $M$.
Define a spinor-valued differential 1-form\/ $\cone{\Theta}_\pm$ on the
$\eps$-metric cone $\cone{M}$ over $M$ by
\begin{align} \label{eq:cone32ksf}
\begin{split}
  \cone{\Theta}_\pm
    & = r \, \cone{\pi}^\Tw (\pi_2^* (F_\pm (\Phi')))
\\
    & = r \big( \pi_2^* (F_\pm (\Phi'))
      + \tfrac{1}{n+1} \, \cclf \wedge (\cclv \intp
          \pi_2^* (F_\pm (\Phi')) \big),
\end{split}
\end{align}
where $\cclf$ and $\cclv$ are the Clifford multiplication form and its
metric dual on $\cone{M}$.
Then\/ $\cone{\Theta}_\pm$ is primitive by construction, and it is
parallel with respect to\/ $\cone{\covdg}$ if and only if\/ $\Phi'$
is a solution of \eqref{eq:32ksf} with the Killing number
$a' = \pm \tfrac{3}{2} \sqrt{\eps}$.

Conversely, any parallel primitive (half\nh) spinor-valued 1-form\/
$\cone{\Theta}$ on $\cone{M}$, in particular,
$\cone{\Theta} \in \Sec(\Spnr^\Tw \cone{M})$ for $n$ even, and\/
$\cone{\Theta} \in \Sec(\Spnr^{\Tw\pm} \cone{M})$ for $n$ odd,
arises this way with\/ $\Phi'_\pm$ given by
\begin{align} \label{eq:cone32ksfinv}
  \Phi'_\pm
    & = \tfrac{1}{2} \big(
        F_\mp (\pi^\Tan (\cone{\Theta} |_M))
        \pm \sqrt{\eps} \, \clf \wedge
            F_\mp (\pi^\Norm (\cone{\Theta} |_M))
      \big).
\end{align}
\end{prop}

\begin{prf}
The proof is based on repeated application of the
relationship between the Clifford multiplication forms $\clf$ and
$\cclf$ on $M$ and $\cone{M}$, respectively.
Taking into account \eqref{eq:cla} and \eqref{eq:coneg}, we easily deduce
\begin{align*}
  \cclf
    & = \dr \otimes ({\vr\clp}) + r \, \pi_2^* (\clf),
  & \cclv
    & = \eps \, \vr \otimes ({\vr\clp}) + \tfrac{1}{r} \, \pi_2^*(\clv).
\end{align*}
We recall that $\pi_2^*$ denotes the pull-back along the projection
$\pi_2 \colon \cone{M} \to M$, and we also canonically identify the
pull-back bundles $\pi_2^* (\cone{\Tan} M)$ and
$\pi_2^* (\cone{\Spnr} M)$ with $\Tan \cone{M}$ and $\Spnr \cone{M}$,
respectively.

First we define the 1-form $\cone{\Theta}_\pm$ on $\cone{M}$ by
\eqref{eq:conesksf} and show that it is primitive if and only if the
forms $\Phi$ and $\Xi$ on $M$ are related by \eqref{eq:32ksfprimit},
\begin{align*}
  \cclv \intp \cone{\Theta}_\pm
    & = \cclv \intp (
        \dr \wedge \pi_2^* (F_\pm (\Phi))
        + r \, \pi_2^* (F_\pm (\Xi))
      )
\\
    & = \vr \intp (
        \dr \wedge \pi_2^* (\eps \, \vr \clp F_\pm (\Phi))
      )
      + \pi_2^* (\clv \intp F_\pm (\Phi))
\\
    & = \pi_2^* (F_\mp (\pm \sqrt{\eps} \, \Phi + \clv \intp \Xi))
\\
    & = \pi_2^* (F_\mp (2 a \Phi + \clv \intp \Xi)),
\end{align*}
where we used \eqref{eq:F}, \eqref{eq:cla} and
$a = \tfrac{1}{3} \, a' = \pm \tfrac{1}{2} \sqrt{\eps}$.
The claimed correspondence between $\Phi'$ and $\cone{\Theta}$ now
follows from Propositions \ref{conesksf} and \ref{32ksfcor}.
The formula \eqref{eq:cone32ksfinv} follows immediately by substituting
\eqref{eq:conesksfinv} for $\Phi$ and $\Xi$ into \eqref{eq:32ksfcor2}.

As for \eqref{eq:cone32ksf}, we substitute \eqref{eq:32ksfcor1} for
$\Phi$ and $\Xi$ into \eqref{eq:conesksf} and compute,
\begin{align*}
\begin{split}
  \cone{\Theta}_\pm
    & =
      \mp \tfrac{1}{\sqrt{\eps} (n+1)} \, \dr \wedge
          \pi_2^* (F_\pm (\clv \intp \Phi'))
      + {}
\\
  & \quad
      + r \, \pi_2^* \big(F_\pm \big(\Phi'
          + \tfrac{1}{n+1} \, \clf \wedge (\clv \intp \Phi')
        \big) \big)
\end{split}
\\
\begin{split}
    & = \tfrac{1}{n+1} \, r \, \dr \wedge (\vr \clp
          (\cclv \intp \pi_2^* (F_\pm (\clv \intp \Phi'))))
      + {}
\\
    & \quad
      + r \, \pi_2^* (F_\pm (\Phi'))
      + \tfrac{1}{n+1} \, r^2 \pi_2^* (\clf) \wedge
          (\cclv \intp \pi_2^* (F_\pm (\Phi')))
\end{split}
\\
    & = r \big(\pi_2^* (F_\pm (\Phi'))
        + \tfrac{1}{n+1} \, \cclf \wedge (\cclv \intp
            \pi_2^* (F_\pm (\Phi')))
      \big),
\end{align*}
which completes the proof.
\end{prf}

Now we return to the example space $\Me$ of constant curvature.
Substituting \eqref{eq:conefproj}, \eqref{eq:F} and \eqref{eq:ccnormal}
into \eqref{eq:cone32ksfinv} we get explicit formulas for Killing
spinor-valued 1-forms $\Phi'_\pm$ on $\Me$ which have the Killing number
$a' = \pm \tfrac{3}{2} \sqrt{\eps}$,
\begin{align} \label{eq:cc32ksf}
\begin{split}
  \Phi'_\pm (x)
    & = \tfrac{1}{2} \big(
        (1 \pm \sqrt{\eps} \, x) \clp (
            \cone{\Theta} - \dr \wedge (x \intp \cone{\Theta})
          )
        \pm {}
\\ & \qquad
        \pm \sqrt{\eps} \, \clf \wedge (
            (1 \pm \sqrt{\eps} \, x) \clp (x \intp \cone{\Theta})
          )
      \big),
\end{split}
\end{align}
where $\cone{\Theta} \in \Spnr^\Tw_{n+1}$ for $n$ even, and
$\cone{\Theta} \in \Spnr^{\Tw\pm}_{n+1}$ for $n$ odd, regarded as a
constant section of $\Spnr^\Tw \cone{M}$ or $\Spnr^{\Tw\pm} \cone{M}$
respectively.
Note that these additional solutions are \emph{not} in the span of
tensor products $\alpha \otimes \Psi$ of a Killing-Yano form and a
Killing spinor.
This is simply due to the fact that there are no nontrivial Killing
spinors on $\Me$ with Killing number
$a' = \pm \tfrac{3}{2} \sqrt{\eps}$.

By Propositions \ref{ccksfint1} and \ref{32ksfcor}, there are no
nontrivial Killing spinor-valued 1-forms on $\Me$ with Killing number
\begin{align*}
  a' \ne 0, \,
    \pm \tfrac{1}{2} \sqrt{\eps}, \,
    \pm \tfrac{3}{2} \sqrt{\eps}
\end{align*}
because, by Proposition \ref{ccsksf0}, there are no nontrivial special
Killing 0-forms on $\Me$ with Killing number
$a \ne \pm \tfrac{1}{2} \sqrt{\eps}$.
Finally we resolve the remaining case $a' = 0$ of the equation
$\eqref{eq:32ksf}$ and hence complete our discussion of all
Killing spinor-valued forms on $\Me$.
In this case the equation simply requires $\covdg$-parallel
spinor-valued 1-forms.

\begin{prop}
There are no nontrivial\/ $\covdg$-parallel spinor-valued 1-forms on the
space $\Me$.
\end{prop}

\begin{prf}
Suppose that $\Phi$ is a $\covdg$-parallel spinor-valued 1-form.
Then the first integrability condition requires that $\Phi$ is
annihilated by the curvature of $\covdg$ and thus we have by
\eqref{eq:cccurv}
\begin{align*}
\begin{split}
  0 = \curvg_{X,Y} \Phi = \eps \rho(X \wedge Y) \, \Phi .
\end{split}
\end{align*}
Hence $\Phi = 0$ since the representation of the spin Lie algebra on
spinor-valued 1\nh forms contains no trivial summands.
Alternatively, we can compute the operator $r_2$ as in Lemma
\ref{ccksfint02},
\begin{align*}
  r_2 (\curvg \Phi)
    & = - \tfrac{1}{2} \, \eps (n - 1) \, \clv \intp \Phi,
\end{align*}
and then the operator $q$ as in Proposition \ref{ccsksf0},
\begin{align*}
  q (\curvg \Phi)
    & = - \big( n - \tfrac{1}{2} \big) \, \Phi
      - \tfrac{1}{2} \, \clf \wedge (\clv \intp \Phi),
\end{align*}
and the claim follows.
\end{prf}

\section{Final remarks and comments}

We have shown that application of the integrability conditions revealed
unexpected Killing spinor-valued $1$-forms on spaces of constant
curvature.
Our results can be regarded as the first example resulting from the investigation
of Killing spinor-valued forms that is not implied by known results on
Killing-Yano forms and Killing spinors.
Apparently, the application towards explicit examples is computationally
rather complicated to do by hand even in the simplest case of spaces of
constant curvature.

However, as the equation \eqref{eq:parint} suggests, all the
integrability conditions can be applied algorithmically, and in
many cases this approach is sufficient to completely determine the
space of solutions.
The second author has implemented an algorithm for solving the three
types of Killing equations on homogeneous spaces using a computer
algebra system, see \cite{zim17}.
In fact, the additional solutions on spaces of constant curvature were
originally discovered this way.
Computed examples include the \emph{Berger spheres} in dimensions
$3$, $5$, $7$ which are Sasakian manifolds, the \emph{Aloff-Wallach space}
$N(1,1)$ which is a nontrivial 3-Sasakian manifold, and the seven-sphere
equipped with $\G_2$-structure.
As a result a new type of solutions appears in the 3-Sasakian case, and
this case will be discussed in a separate article.

It is worth of notice that the relationship between the existence of
Killing spinor-valued forms and the \emph{Einstein manifolds} is not
clear.
Contrary to Killing spinors in the Riemannian case, the Einstein
condition imposed on curvature is not a direct consequence of the
integrability conditions.
On the other hand, there are not known counterexamples.
For example, computer aided computations produced no solutions on
Berger spheres with non-Einstein metrics.